\documentclass[11pt,reqno]{amsart}
\usepackage{amssymb}
\usepackage{amscd}
\usepackage{amsmath}

\theoremstyle{plain}
\newtheorem{theorem}{Theorem}[section]
\newtheorem{lemma}[theorem]{Lemma}
\newtheorem{corollary}[theorem]{Corollary}
\newtheorem{proposition}[theorem]{Proposition}

\theoremstyle{definition}
\newtheorem{remark}[theorem]{Remark}
\newtheorem{remarks}[theorem]{Remarks}
\newtheorem{definition}[theorem]{Definition}
\newtheorem{example}[theorem]{Example}
\newtheorem{examples}[theorem]{Examples}

\newcommand{\bq}{\mathbf Q}
\newcommand{\bz}{\mathbf Z}
\newcommand{\br}{\mathbf R}
\newcommand{\bc}{\mathbf C}

\newcommand{\F}{\mathcal F}

\newcommand{\calP}{\mathcal P}
\newcommand{\calC}{\mathcal C}
\newcommand{\calO}{\mathcal O}
\newcommand{\calL}{\mathcal L}
\newcommand{\A}{\mathfrak A}
\newcommand{\fA}{\mathfrak A}
\newcommand{\B}{\mathfrak B}

\newcommand{\frP}{\mathfrak P}

\newcommand{\Is}{\mathrm{Is}}
\newcommand{\nr}{\mathrm{nrd}}

\newcommand{\lra}{\longrightarrow}

\newcommand{\bbA}{{\mathbb A}}
\newcommand{\cA}{{\mathcal A}}

\newcommand{\bC}{{\mathbf C}}
\newcommand{\fC}{{\mathfrak C}}

\newcommand{\bG}{{\mathbf G}}

\newcommand{\frp}{{\mathfrak p}}

\newcommand{\hqed}{{\hfill \qed}}
\newcommand{\bH}{{\mathbf H}}
\newcommand{\cH}{{\mathcal H}}

\renewcommand{\Im}{{\text{Im}}}
\newcommand{\bj}{{\mathbf j}}
\newcommand{\vk}{{\varkappa}}

\newcommand{\opsi}{{\boldsymbol{\psi}}} 

\newcommand{\oK}{{\overline{K}}}
\renewcommand{\L}{{\mathcal L}}
\newcommand{\oL}{\overline{{\mathcal L}}}
\newcommand{\M}{{\mathfrak M}}
\newcommand{\cM}{{\mathcal M}}

\newcommand{\cp}{\mathcal P}
\renewcommand{\O}{{\mathcal O}}

\newcommand{\cP}{{\mathcal P}}

\newcommand{\Q}{{\mathbf Q}}
\newcommand{\bQ}{{\mathbf Q}}

\renewcommand{\r}{{\mathbf r}}
\newcommand{\bR}{{\mathbf R}}

\renewcommand{\top}{{\mathrm{top}}}

\newcommand{\Z}{{\mathbf Z}}
\newcommand{\bZ}{{\mathbf Z}}

\DeclareMathOperator{\AC}{AC}

\DeclareMathOperator{\Aut}{Aut}

\DeclareMathOperator{\codim}{codim}
\DeclareMathOperator{\Coker}{Coker}
\DeclareMathOperator{\Cok}{Cok}

\DeclareMathOperator{\Det}{Det} 

\DeclareMathOperator{\End}{End} 

\DeclareMathOperator{\Ext}{Ext} 
\DeclareMathOperator{\fin}{fin}
\DeclareMathOperator{\Gal}{Gal} 
\DeclareMathOperator{\GL}{GL}
\DeclareMathOperator{\Hom}{Hom} 
\DeclareMathOperator{\Image}{Im}
\DeclareMathOperator{\Ker}{Ker} 

\DeclareMathOperator{\Map}{Map} 

\DeclareMathOperator{\Pic}{Pic} 
\DeclareMathOperator{\PPic}{PPic}
\DeclareMathOperator{\hPic}{\widehat{Pic}}
\DeclareMathOperator{\rag}{rag}

\DeclareMathOperator{\Spec}{Spec}

\DeclareMathOperator{\Ind}{Ind}

\DeclareMathOperator{\im}{im}

\DeclareMathOperator{\id}{id}

\title{On Twisted Forms and Relative Algebraic $K$-theory}
\author{A. Agboola \and D. Burns}
\address{Dept. of Math., University of California at Santa
Barbara, CA 93106, U.S.A.}
\email{agboola@math.ucsb.edu}
\address{Dept. of Math., King's College London, Strand,
London WC2R 2LS, England} 
\email{david.burns@kcl.ac.uk}
\date{Final version, March 23, 2005}
\begin{document}

\begin{abstract}
This paper introduces a new approach to the study of certain aspects
of Galois module theory by combining ideas arising from the study of
the Galois structure of torsors of finite group schemes with
techniques coming from relative algebraic $K$-theory.
\end{abstract}

\maketitle

\newcommand{\Ze}{\mathbf{Z}}
\newcommand{\Ce}{\mathbf{C}}
\renewcommand{\Re}{\mathbf{R}}
\renewcommand{\Im}{ \mathrm{Im} }

\newcommand{\CG}{ {\Ce [G]} }
\newcommand{\ZG}{ {\Ze [G]} }
\newcommand{\CN}{ {\calC_N} }
\newcommand{\DS}{ {\Delta S} }
\newcommand{\cok}{\text{cok}}
\newcommand{\diff}{\frd_k}
\newcommand{\indGp}{\text{ind}_{G_\frP}^G}
\newcommand{\lmt}{\longmapsto}
\newcommand{\ra}{\rightarrow}
\newcommand{\nI}{ {\frn_I} }
\newcommand{\nJ}{ {\frn_J} }
\newcommand{\Ol}{{\calO_l}} 
\newcommand{\QG}{ {\bq G} }
\newcommand{\RG}{ {\Re [G]} }
\newcommand{\QpG}{ {\bq_p[G]} }
\newcommand{\QlG}{ {\bq_l G} }
\newcommand{\Zp}{ {\Ze_p} }
\newcommand{\ZpG}{ {\Ze_p[G]} }
\newcommand{\tensor}{ {\otimes} }
\newcommand{\tensorZ}{ {\otimes_\Ze} }
\newcommand{\tensorQ}{ {\otimes_\bq} }
\newcommand{\ZS}{\Ze S}
\newcommand{\Zl}{ {\Ze_l} }
\newcommand{\ZlG}{ {\Ze_l[G]} }
\newcommand{\TLK}{ T_{L/K} }
\newcommand{\WLKinf}{ w_\infty(L/K) }
\newcommand{\HL}{ {H_L} }
\newcommand{\HLZ}{ {H_{L, \Ze}}}
\newcommand{\HLKZ}{ {H_{L/K, \Ze}}}
\newcommand{\SL}{ {\Sigma(L)} }
\newcommand{\SK}{ {\Sigma(K)} }
\newcommand{\piLK}{ {\pi_{L/K} } }
\newcommand{\tauLK}{ {\tau_{L/K} } }

\maketitle

\section{Introduction}

Let $R$ be a Dedekind domain with field of fractions $F$, and suppose
that $G \to \Spec(R)$ is a finite, flat, commutative group scheme. In
recent years, there has been some interest in the study of Galois
structure invariants attached to torsors of $G$. The original
motivation for this work arose from the study of the Galois module
structure of rings of integers: if $X \to \Spec(R)$ is a $G$-torsor,
then we may write $X = \Spec(\fC)$, and in many cases, the algebra
$\fC$ may be viewed as being an order in the ring of integers of some
(in general wildly ramified) field extension of $F$. The case in which
$G$ is a torsion subgroup scheme of an abelian variety is particularly
interesting: here the corresponding torsors are obtained by dividing
points in the Mordell-Weil group of the abelian variety, and
this enables one to relate questions concerning the Galois structure
of rings of integers to the arithmetic of abelian varieties. This
approach was first introduced by M. J. Taylor (see \cite{T1},
\cite{T2}, \cite{BT}), and has since been developed in greater
generality by several authors (see for example \cite{A3}, \cite{AP},
\cite{By}, \cite{P1}).

The main goal of this paper is to show that combining ideas arising
from the study of Galois structure invariants attached to torsors of
$G$ with techniques from relative algebraic $K$-theory yields a
natural refinement and reinterpretation of several different aspects
of Galois module theory. At the same time, we shall also see that this
approach also gives fresh insight into a number of old results. We
remark that techniques involving relative algebraic $K$-theory have
already played an essential role in the formulation of an
equivariant version of the Tamagawa Number Conjecture of Bloch and
Kato, and in the description of its consequences for Galois module
theory (see \cite{Bu1}, \cite{BuF}).

We now give an outline of the contents of this paper. In
\S\ref{S:twisted}, we recall the definition of a fibre product
category, and we introduce the notion of a {\it generalised twisted
form}. This is a generalisation of the usual notion of a twisted form
or principal homogeneous space of a Hopf algebra (see Remark
\ref{R:motivation}). An example of such an object may be given as
follows. Suppose that $\A$ is a finitely generated $R$-algebra which
spans a semisimple $F$-algebra $A$, and let $E$ be an extension of
$F$. Then a generalised $E$-twisted $\A$-form is given by a triple
$[M,N;\iota]$, where $M$ and $N$ are finitely generated, locally free
$\A$-modules, and $\iota: M \otimes_R E \xrightarrow{\sim} N \otimes_R
E$ is an isomorphism of $A \otimes_F E$-modules. We discuss two
particular examples of generalised twisted forms that provide
motivation for our later work. The first example arises via finite
Galois extensions of number fields, while the second is constructed
using torsors of finite, flat, commutative group schemes.

It follows from the definition that generalised twisted forms may be
classified by appropriate relative algebraic $K$-groups. (For example,
the generalised twisted forms $[M,N;\iota]$ described above yield
classes in a relative algebraic $K$-group that we denote by
$K_0(\A,E)$.) In \S\ref{S:ktheory}, we recall the basic facts
concerning relative algebraic $K$-theory that we use, and we then give
a Fr\"ohlich-type `Hom-description' of relative algebraic $K$-groups
of the form $K_0(\A,E)$. This description may itself be of some
independent interest.

Suppose now that $F$ is a number field. In \S\ref{S:metrised}, we
define the notion of a metric $\mu$ on a finitely generated, locally
free $\A$-module $M$. Following \cite{CPT}, we then define an
arithmetic classgroup $\mathrm{AC}(\A)$ which classifies metrised
$\A$-modules $(M,\mu)$. There is a natural homomorphism $\partial_{\A,
F^c}$ from the relative algebraic $K$-group $K_0(\A,F^c)$ to
$\mathrm{AC}(\A)$ (where $F^c$ denotes an algebraic closure of
$F$). Proposition \ref{pr} implies that if $(M,N;\iota_1)$ is a
generalised $F^c$-twisted $\A$-form with $M$ and $N$ both locally free
$\A$-modules, then for \textit{any} metric $\mu$ on $N$, the image of
the class $(M,N;\iota_1)$ under $\partial_{\A, F^c}$ measures the
difference between the elements of $\mathrm{AC}(\A)$ which are
associated to $(M,\iota_{1}^{*}(\mu))$ and $(N,\mu)$, where
$\iota_{1}^{*}(\mu)$ is the metric on $M$ which is obtained from $\mu$
via pullback by $\iota_1$. This result is of interest since, in many
of the examples that occur in Galois module theory, canonical metrics
on locally free $\A$-modules can be shown to be equal to, for
instance, the pullback of the trivial metric on $\A$ by a natural
homomorphism.

In \S\ref{S:classical}, we give two applications of the approach
described in this paper to the study of Galois structure of rings of
integers in tamely ramified extensions. For the first application, we
combine our work with results of Bley and the second-named author from
\cite{BB} and we obtain a natural strengthening of a result of
Chinburg, Pappas and Taylor concerning the Hermitian Galois structure
of rings of integers. For the second application, we describe how
certain natural torsion Galois module invariants first introduced by S. Chase
(see \cite{Ch}) are related to certain `equivariant discriminants'
that take values in a relative algebraic $K$-group.

The remainder of the paper is devoted to describing how the methods we
develop may be applied to the study of the Galois structure of torsors
in several different contexts. In \S\ref{S:torsors}, we discuss
reduced resolvends associated to torsors. This notion was first
introduced by L. McCulloh in a different setting (see \cite{Mc}). We
use a different approach from McCulloh's (see \cite{A3}), and we give
a new characterisation of reduced resolvends as being primitive
elements of a certain algebraic group.

We then introduce a natural refinement of the class invariant
homomorphism first studied by W. Waterhouse in \cite{W}. Our refined
homomorphism takes values in a suitable relative algebraic
$K$-group. We show that the homomorphism is injective, and that its
image admits a precise functorial description in almost all cases of
interest. This extends a result of the first-named author in
\cite{A2}.

Finally, we explain how the approach described in this paper enables
one to refine McCulloh's results on realisable classes of rings of
integers of tame extensions by considering invariants that lie an
appropriate relative algebraic $K$-group, rather than in a locally
free classgroup. This gives a relationship between McCulloh's
realisability results, and the work of Chase concerning certain
torsion Galois modules. It also yields realisability results in
certain arithmetic classgroups. We find that, in general, the
collection of realisable classes in the relative algebraic $K$-group
does not form a group.
\medskip

\noindent{}{\bf Acknowledgements.} We would like to thank the
Mathematics Department of Harvard University for its hospitality while
a part of this work was carried out. We are also very grateful to the
referee of this paper for a number of very helpful remarks. The
research of the first-named author was partially supported by NSF
grants DMS-9700937 and DMS-0070449.
\medskip

\noindent{}{\bf Notation.}
Throughout this paper, all modules are left modules, unless explicitly
stated otherwise. 

For any field $L$, we write $L^c$ for an algebraic closure of $L$, and
we set $\Omega_L:= \Gal(L^c/L)$. If $L$ is either a number field or a
local field, then $O_L$ will denote its ring of integers.

If $L$ is a number field, we write $S_f(L)$ (resp. $S_{\infty}(L)$)
for the set of finite (resp. infinite) places of $L$. If $v$ is any
place of $L$, then we write $L_v$ for the local completion of $L$ at
$v$. For any $O_L$-module $P$, we shall usually set $P_v:= P
\otimes_{O_L} O_{L_v}$.

If $S_1 \subset S_2$ are rings, and if $A$ is any $S_1$-algebra, then
we often write $A_{S_2}$ for $A \otimes_{S_1} S_2$. We also often use
similar notation $M_{S_2}:=M \otimes_{S_1} S_2$ for any $S_1$-module $M$.

The symbol $\zeta(A)$ denotes the centre of $A$. 

For any isomorphic $A$-modules $M$ and $N$, we write $\Is_A(M,N)$ for
the set of $A$-equivariant isomorphisms $M \xrightarrow{\sim}
N$. 

Throughout this paper, $R$ denotes a Dedekind domain with field of
fractions $F$. We write $R^c$ for the integral closure of $R$ in
$F^c$. We also allow the possibility that $R=F$, in which case of
course $R^c=F^c$. 

Finally, we remark that if $G$ and $H$ are two groups, then we
sometimes use the notation $g \times h$ (rather than $(g,h)$) for an
element of the product $G \times H$. \hqed

\section{Twisted Forms and Fibre Product Categories} \label{S:twisted}

\medskip
\subsection{Fibre product categories}\label{fpci} 
Let $F_i:\cp_i\rightarrow\cp_3$, $i\in \{1,2\}$ be functors between
categories and consider the {\em fibre product category}
$\cp_4:=\cp_1\times_{\cp_3}\cp_2$. Explicitly, $\cp_4$ is the
category with objects $(L_1,L_2;\lambda)$ where $L_i$ is an object of
$\cp_i$ for $i \in \{1,2\}$ and
$\lambda:F_1(L_1)\xrightarrow{\sim}F_2(L_2)$ is an isomorphism in
$\cp_3$. Morphisms $\alpha:(L_1,L_2;\lambda) \rightarrow
(L_1',L_2';\lambda')$ in $\cp_4$ are pairs
$\alpha=(\alpha_1,\alpha_2)$ with $\alpha_i\in\Hom_{\cp_i}(L_i,L_i')$
so that the diagram
\[
\begin{CD}
F_1(L_1) @>F_1(\alpha_1)>> F_1(L_1')\\
@VV\lambda V @VV\lambda'V\\ F_2(L_2) @>F_2(\alpha_2)>> F_2(L_2')
\end{CD}
\]
in $\cp_3$ commutes. Such a morphism $\alpha$ is an isomorphism in
$\cp_4$ if and only if the morphisms $\alpha_1$ and $\alpha_2$ are
both isomorphisms. (The reader may consult \cite[Chapter VII,
\S3]{bass} for more details concerning fibre product categories.)

\subsection{Generalised twisted forms}

For any ring $\Lambda$ we write $\calP_\Lambda$ for the category
of finitely generated projective $\Lambda$-modules.

Recall that $R$ is a Dedekind domain, with fraction field $F$. Suppose
that $\A$ is a finitely generated $R$-algebra which spans a semisimple
$F$-algebra $A$. For any extension $\Lambda$ of $R$ we write
$\calP_\A\times_\Lambda\calP_\A$ for the fibre product category with
$F_1$ and $F_2$ both equal to the scalar extension functor
$-\otimes_R\Lambda$ from $\calP_\A$ to $\calP_{\A\otimes_R\Lambda}$.

\begin{definition} \label{D:twisted}
We shall refer to an object of the category $\calP_{\A} \times_\Lambda
\calP_{\A}$ as a {\it generalised $\Lambda$-twisted $\A$-form}. We
abbreviate this terminology to {\it generalised twisted form}, or even
just {\it twisted form}, when both $\A$ and $\Lambda$ are clear from
the context. \hqed
\end{definition}

\begin{example}\label{class} (Classical) 
Suppose that $F$ is a number field, and let $L/F$ be any finite
extension. We write $N$ for the normal closure of $L$ over $F$ and we
let $\Sigma_F(L,N)$ denote the set of distinct $F$-embeddings of $L$
into $N$.  Set $H_F^L := \prod_{\Sigma_F(L,N)}\Ze$, and write
$\Aut(L/F)$ for the group of $F$-automorphisms of $L$. Then there is a
canonical $N[\Aut(L/F)]$-equivariant isomorphism
\[ 
\pi^{L,N}_F: L\otimes_F N \xrightarrow{\sim} H_F^L\otimes_\bz N
\]
which is given by mapping each primitive tensor $l\otimes n$ to the
vector $(\sigma(l) n)_\sigma$.

Now suppose that $S$ is any subring of $F$ and let $G$ be a subgroup
of $\Aut(L/F)$.  Then for any projective $S[G]$-submodule $\calL$ of $L$ such
that $\calL \otimes _S F = L$, the triple $(\calL,H_F^L\otimes_\bz
S;\pi^{L,N}_F)$ is a generalised
$N$-twisted $S[G]$-form.

In special cases there are canonical choices of $\calL$. For example,
if $S=F$, then one can take $\calL = L$. In addition, if $L$ is a
tamely ramified Galois extension of a field $K$ with $\Gal(L/K)=G$ and
$S = O_F$ for any subfield $F$ of $K$, then one can take $\calL$ to be
any $G$-stable fractional ideal of $O_L$ (see \cite{U}). \hqed
\end{example}

\begin{example}\label{geom} (Geometrical) 
Let $R$ be any Dedekind domain, with field of fractions $F$. (We allow
$R=F$.) Suppose that $G \to \Spec(R)$ is a finite, flat, commutative
group scheme, and let $G^*$ denote its Cartier dual. Write $\A:=
\calO_{G^*}$ and $A:= \A \otimes_R F$. It is shown by Waterhouse in
\cite[Sections 1 and 2]{W} (see also \cite[\'expos\'e VII]{SGA7} or
\cite{P}) that there is a canonical isomorphism of groups
\begin{equation*}
H^1(\Spec(R),G) \simeq \Ext^1(G^*,\bG_m).
\end{equation*}
This implies that given any $G$-torsor $\pi:X \to \Spec(R)$, we can
associate to it a canonical commutative extension
\begin{equation} \label{E:extension}
1 \to \bG_m \to G(\pi) \to G^* \to 1.
\end{equation}
The scheme $G(\pi)$ is a $\bG_m$-torsor over $G^*$, and we write
$\L_{\pi}$ for its associated $G^*$-bundle.

Let $\pi_0: G \to \Spec(R)$ denote the trivial $G$-torsor. Then, over
$\Spec(R^c)$, the $G$-torsors $\pi$ and $\pi_0$ become isomorphic,
i.e. there is an isomorphism $X \otimes_R R^c \xrightarrow{\sim} G
\otimes_R R^c$ of schemes with $G$-action. Hence, via the
functoriality of the construction in \cite[Section 2]{W}, we obtain an
$\A \otimes_R R^c$-equivariant isomorphism
\begin{equation} \label{E:splitting}
\xi_{\pi}: \L_{\pi} \otimes_R R^c \xrightarrow{\sim} \A \otimes_R R^c.
\end{equation}
We refer to $\xi_{\pi}$ as a {\it splitting isomorphism} associated to
$\pi$.  Then the triple $(\L_\pi,\A;\xi_\pi)$ is a generalised
$R^c$-twisted $\A$-form. \hqed
\end{example}

\begin{remark} \label{R:motivation}
The line bundle $\L_{\pi}$ in Example \ref{geom} may be described
explicitly in the following way (see \cite[Section 4]{W}). For any
$G$-torsor $\pi:X \to \Spec(R)$ as above, the structure sheaf $\O_X$
is an $\O_G$-comodule, and so it is also an $\A$-module (see
\cite[Proposition 1.3]{CEPT}). As an $\A$-module, $\O_X$ is locally
free of rank one, and so it gives a line bundle $\cM_{\pi}$ on
$G^*$. Then it may be shown that $\L_{\pi} = \cM_{\pi} \otimes
\cM_{\pi_0}^{-1}$. The line bundle $\cM_\pi$ is a principal
homogeneous space, or twisted form of the Hopf algebra $\A$ in the
sense of, for example \cite{T1} or \cite{BT}, and this is the
motivation for Definition \ref{D:twisted} above. \hqed
\end{remark}

\section{Relative algebraic $K$-theory} \label{S:ktheory}

\subsection{Basic theory} 
Let $\Lambda$ be any extension of $R$, and let $\A$ be a finitely
generated $R$-algebra which spans an $F$-algebra $A$. In this
subsection, we shall recall some basic facts concerning the $K$-theory
of categories of the form $\calP_{\A} \times_{\Lambda}
\calP_{\A}$. Further details of these results may be found in
\cite[Chapter 15]{Sw} and \cite[Section 40B]{CR}.

\begin{definition} \label{D:ses}
A short exact sequence of objects in $\calP_{\A} \times_{\Lambda}
\calP_{\A}$ is a sequence
\begin{equation} \label{E:sequence}
0 \to (L_1,L_2;\lambda) \xrightarrow{(\alpha_1,\alpha_2)}
(L'_{1},L'_{2};\lambda') \xrightarrow{(\alpha'_{1},\alpha'_{2})}
(L''_{1},L''_{2};\lambda'') \to 0
\end{equation}
such that $(\alpha_1,\alpha_2)$ and $(\alpha'_{1},\alpha'_{2})$ are
morphisms in $\calP_{\A} \times_{\Lambda} \calP_{\A}$, and such that
\begin{equation*}
0 \to L_1 \xrightarrow{\alpha_1} L'_{1}\xrightarrow{\alpha'_{1}}
L''_{1} \to 0 \quad \text{and} \quad
0 \to L_2 \xrightarrow{\alpha_2} L'_{1}\xrightarrow{\alpha'_{2}}
L''_{2} \to 0
\end{equation*}
are exact sequences of $R$-modules. \hqed
\end{definition}

\begin{definition} \label{D:relk}
We write $[L_1,L_2;\lambda]$ for the isomorphism class of an object
$(L_1,L_2;\lambda)$ of $\calP_{\A} \times_{\Lambda} \calP_{\A}$. We
define the relative algebraic $K$-group $K_0(\A,\Lambda)$ to be the
abelian group with generators the isomorphism classes of objects of
$\calP_{\A} \times_{\Lambda} \calP_{\A}$ and relations
 
(i) $[L'_{1},L'_{2};\lambda'] = [L_1,L_2;\lambda] +
[L''_{1},L''_{2};\lambda'']$ for each exact sequence (3);
 
(ii) $[M_1,M_2; \eta\circ\eta'] = [M_1,M_3; \eta'] +
[M_3,M_2;\eta]$.
 
It may be shown that every element of $K_0(\A,\Lambda)$ is of the form
$[L_1,L_2 ;\lambda]$ (see \cite[Lemma 15.6]{Sw}). However, the natural
map from the set of isomorphism classes of objects of $\calP_{\A}
\times_{\Lambda} \calP_{\A}$ to $K_0(\A,\Lambda)$ is very far from
being injective (for example, the image of the element
$[\A^n,\A^n;1_{\A_{\Lambda}}]$ under this map is zero for every
positive integer $n$). \hqed
\end{definition}

For any ring $\Lambda$ we write $K_0T(\Lambda)$ for the
Grothendieck group of the category of $\Lambda$-modules which are
both finite and of finite projective dimension. Then there are natural
isomorphisms
\begin{equation}\label{rt} 
K_0(\A,F) \xrightarrow{\sim}
K_0T(\A) \xrightarrow{\sim} \bigoplus_v K_0T(\A_v)
\end{equation}
where $v$ runs over all finite places of $R$ (see the discussion
following (49.12) in \cite{CR}). The first isomorphism in \eqref{rt}
is defined in the following way. If $M$ is any $\fA$-module which is
both finite and of finite projective dimension, then (since $R$ is a
Dedekind domain), there exists an exact sequence of $\fA$-modules of
the form
$$
0 \to O \xrightarrow{\varphi} Q \to M \to 0
$$
in which $P$ and $Q$ are both finitely generated and projective. The
first isomorphism of \eqref{rt} sends $(P,\varphi,Q)$ to the class of
$M$. 

For any ring $S$ we write $K_i(S)$ ($i=0,1$) for the algebraic
$K$-group in degree $i$ of the exact category $\cP_S$. We recall that
if $R \to \Lambda$ is any homomorphism of rings, then there is a long
exact sequence of relative algebraic $K$-theory
\begin{equation} \label{krel} 
K_1(\A) \xrightarrow{\partial^{2}_{\fA,\Lambda}}  K_1(\A_{\Lambda})
\xrightarrow{\partial^{1}_{\A,\Lambda}}   
K_0(\A,\Lambda) \xrightarrow{\partial^{0}_{\A,\Lambda}}  
K_0(\A) \to K_0(\A_{\Lambda})
\end{equation} 
(see \cite[Theorem 15.5]{Sw}). Here the homomorphism
$\partial^{2}_{\fA,\Lambda}$ is the natural scalar extension morphism. The
homomorphism $\partial^{0}_{\A,\Lambda}$ is defined by
\begin{equation*}
\partial^{0}_{\A,\Lambda} ([L_1,L_2;\lambda]) = [L_1] - [L_2].
\end{equation*}
In order to describe $\partial^{1}_{\A,\Lambda}$, we first recall that
$K_1(A_{\Lambda})$ is generated by elements of the form $(V,\phi)$, where $V$
is a finitely generated free $\A_{\Lambda}$-module and $\phi \in
\Is_{\A_{\Lambda}}(V,V)$. If $T$ is any projective $\A$-submodule of $V$ such
that $T\otimes_R \Lambda = V$, then $\partial^{1}_{\A,\Lambda}$ is defined by
setting
\[ 
\partial^{1}_{\A,\Lambda}((V,\phi)) = [T,T;\phi].
\]
This definition is independent of the choice of $T$.

\begin{remark}\label{ct}
Let $E$ be any extension of $F$. When $A_E$ is a semisimple algebra,
it is often convenient to compute in $K_1(A_E)$ by means of the
injective `reduced norm' map
\begin{equation*}
 \mathrm{nrd}_{A_E}: K_1(A_E) \to \zeta(A_E)^{\times}.
\end{equation*}
This map sends the element $(V,\phi)$ to the reduced norm of
 $\phi$, considered as an element of the semisimple $E$-algebra
$\End_{A_E}(V)$.

If $\A$ is commutative, then the determinant functors over $A_E$ and
$\A$ combine to induce canonical isomorphisms of $K_1(A_E)/\im
(\partial^{2}_{\fA,E})$ with $A_E^\times /\A^\times$ and of $\ker
[K_0(\A) \lra K_0(A_E)]$ with $\Pic(\A)$. Hence the exact sequence
\eqref{krel} implies that in this case, the group $K_0(\A,E)$ may be
identified with the multiplicative group of invertible $\A$-modules in
$A_E$. \hqed
\end{remark}

\begin{proposition} \label{T:relk}
If $E$ is any extension of $F$, then the following sequence is exact;
\begin{equation} \label{E:relk}
0 \to K_1(A) \xrightarrow{f_1} K_1(A_E) \oplus K_0(\fA,F)
\xrightarrow{f_2} K_0(\fA,E) \to 0.
\end{equation}

Here $f_1(x) = (\partial^{2}_{A,E}(x), -\partial^{1}_{\fA,F}(x))$, and
$f_2(y_1,y_2) = \partial^{1}_{\fA,E}(y_1) + \iota_{\fA,E}(y_2)$, where
$\iota_{\fA,E}: K_0(\fA,F) \to K_0(\fA,E)$ is the natural scalar
extension morphism.
\end{proposition}

\begin{proof} We first observe that $f_1$ is injective because the
scalar extension map $\partial^{2}_{A,E}$ is injective. Next, we note
that upon comparing the exact sequences (\ref{krel}) with $\Lambda=E$
and $\Lambda=F$, one obtains a commutative diagram:
\begin{equation*} \label{E:kdiag}
\minCDarrowwidth1em
\begin{CD} 
K_1(\A) @>{f_3}>> K_1(A_E) @>{\partial^{1}_{\fA,E}}>> K_0(\A ,E )
@>{f_4}>> K_0(\A) 
@>{f_5}>> K_0(A_E)\\  
@\vert @A{\partial^{2}_{A,E}}AA @A\iota_{\fA,E}AA @\vert @A\iota_{A,E}AA\\ 
K_1(\A) @>{f_6}>> K_1(A) @>{\partial^{1}_{\fA,F}}>> K_0(\A ,F) @>{f_7}>>
K_0(\A) @>>> K_0(A). 
\end{CD}
\end{equation*}
In this diagram, the morphisms $\iota_{\fA,E}$ and $\iota_{A,E}$ are
the natural scalar extension morphisms, and are therefore
injective. (We remark that in fact the injectivity of $\iota_{\fA,E}$
follows from the injectivity of $\partial^{2}_{A,E}$, as may be shown
by an easy diagram-chasing argument.)

We now show that $f_2$ is surjective. Suppose that $z \in
K_0(\fA,E)$. Then, as $f_5(f_4(z))=0$ and $\iota_{A,E}$ is injective,
we may choose $z_1 \in K_0(\fA,F)$ such that $f_7(z_1) = f_4(z)$. We
then have that $z - \iota_{\fA,E}(z_1) = \partial^{1}_{\fA,E}(z_2)$
for some $z_2 \in K_1(A_E)$. It is easy to check that $f_2(z_2,z_1) =
z$, and so it follows that $f_2$ is surjective.

It now remains to show that $\Ker(f_2) = \Im(f_1)$. Suppose that $y_1
\in K_1(A_E)$ and $y_2 \in K_0(\fA,F)$ are such that
$$
\partial^{1}_{\fA,E}(y_1) = \iota_{\fA,E}(y_2) = y \in K_0(\fA,E),
$$
say. We have to show that there exists $y' \in K_1(A)$ such that
$\partial^{2}_{A,E}(y') = y_1$ and $\partial^{1}_{\fA,F}(y') = y_2$.

Now $f_5(y) = f_4(\partial^{1}_{\fA,E}(y_1)) = 0$, and so it follows
that $f_7(y_2) = 0$. Hence there exists $y_3 \in K_1(A)$ such that
$\partial^{1}_{\fA,F}(y_3) = y_2$. We now deduce that $f_3(y_4) =
\partial^{2}_{A,E}(y_3) - y_1$ for some $y_4 \in K_1(\fA)$. Set $y' = y_3 -
f_6(y_4)$. It is easy to show that $\partial^{2}_{A,E}(y') = y_1$ and
$\partial^{1}_{\fA,F}(y') = y_2$. This completes the proof.
\end{proof}

\subsection{Hom-descriptions}

Recall that $R$ is a Dedekind domain with field of fractions
$F$. Suppose now that $F$ is a number field, and let $\Gamma$ be a
finite group upon which $\Omega_F$ acts (possibly trivially).  Let
$\A$ be a finitely generated $R$-subalgebra of $F^c[\Gamma]$ such that
there is an equality
\begin{equation}\label{hd}
\A\otimes_R F^c = F^c[\Gamma].
\end{equation}

Under this condition, for any extension $E$ of $F$, we shall give a
description of $K_0(\A,E)$ in terms of idelic-valued functions on the
ring $R_\Gamma$ of $F^c$-valued characters of $\Gamma$. This
description is modeled on the `Hom-descriptions' of class groups
introduced by Fr\"ohlich (cf.  for example \cite[Chapter II]{Fr0}). It
will be useful in later sections, and is itself perhaps also of some
independent interest.
\smallskip

Let $\hat \Gamma$ denote the set of $F^c$-valued characters of
$\Gamma$. The group $\Omega_F$ acts on $R_\Gamma$ according to the
following rule: if $\chi \in \hat \Gamma$ and $\omega \in \Omega_F$,
then for each $\gamma \in \Gamma$ one has $(\omega \circ \chi)(\gamma)
= \omega (\chi(\omega ^{-1}(\gamma)))$.

We fix an embedding of $F^c$ into $E^c$ and we view each element of
$\hat \Gamma$ as taking values in $E^c$. For each element $a$ of
$\GL_n(A_E)$ we define an element $\Det (a)$ of
$\Hom(R_\Gamma,(E^c)^\times)$ in the following way: if $T$ is a
representation over $F^c$ which has character $\phi$, then
\[ 
\Det (a)(\phi) := \det( T(a)).
\]
This definition depends only upon the character $\phi$, and not upon
the choice of representation $T$.

We write $J_f(F^c)$ for the group of finite ideles of the field $F^c$,
and we view $F^\times$ as being a subgroup of $J_f(F^c)$ via the
natural diagonal embedding. If $a$ is any element of $\GL_n(A\otimes_F
J_f(E))$, then the above approach allows one to define an element
$\Det (a)$ of $\Hom(R_\Gamma, J_f(F^c))$. We set 
\begin{equation*}
U_f(\A) := \prod_{v
\in S_f(F)}\A_v^\times \subset (A\otimes_FJ_f(F))^\times.
\end{equation*}
We then define a homomorphism
\begin{equation} \label{E:deltadef}
\Delta_{\A ,E}:
\Det (A^\times) \to \frac{\Hom(R_\Gamma
 ,J_f(F^c))^{\Omega _F} }{\Det(U_f(\A ))} \times
\Det(A_E^\times); \qquad \theta \mapsto (\theta,
\theta^{-1}).
\end{equation}

\begin{theorem}\label{Hom} 
Under the above conditions, there is a natural isomorphism 
\begin{equation} \label{E:homiso} 
h_{\A ,E}: K_0(\A ,E) \xrightarrow{\sim} \mathrm{Cok}(\Delta_{\A
,E}).
\end{equation}
\end{theorem}

\begin{proof} From (\ref{rt}) and the isomorphisms of \cite[Chapter
II, (2.3)]{Fr0}, it follows that there is an isomorphism
$$
h_1: K_0(\fA,F) \xrightarrow{\sim}
\frac{\Hom(R_\Gamma,J_f(F^c))^{\Omega_F}}{\Det(U_f(\fA))}.
$$
Also, the isomorphisms $K_1(A_\Lambda) \cong \Det (A_\Lambda^\times )$
of [loc. cit., Chapter II, Lemma 1.2 and Lemma 1.6] for $\Lambda \in
\{ E,F\}$ induce an isomorphism
$$
h_2: \frac{K_1(A_E)}{\Im(\partial^{2}_{\fA,E})} \xrightarrow{\sim}
\frac{\Det(A_E^\times)}{\Det(A^\times)}.
$$
It follows that there is a natural isomorphism
$$
\Coker(f_1) \simeq \Coker(\Delta_{\fA,E}),
$$
where $f_1$ is defined in the statement of Theorem \ref{T:relk}. The
desired result now follows from the fact that Theorem \ref{T:relk}
implies that there is a natural isomorphism
$$
\Coker(f_1) \xrightarrow{\sim} K_0(\fA,E).
$$
\end{proof}

\begin{remark} An explicit description of the isomorphism $h_{\A ,E}$
may be given as follows. Suppose then that $[X,Y;\xi]$ is an element
of $K_0(\A,E)$ for which $X_F$ and $Y_F$ are free $A$-modules of the
same rank (by [28, Lemma 15.6] any element of $K_0(\A,E)$ is of this
shape). Then for any choice of isomorphism $\theta \in \Is_A(X_F,Y_F)$
one has
\begin{align*} 
[X,Y;\xi] &=[X,Y;\theta_E] + [Y,Y;\xi\circ \theta_E^{-1}] \\ 
&= [X,Y;\theta_E]+ \partial_{\A,E}^1((Y_E,\xi\circ\theta_E^{-1}))
\end{align*}
in $K_0(\A,E)$. The element $\,h_{\A ,E}([X,Y;\xi])$ of
$\Cok(\Delta_{\A,E})$ is then represented by the pair \break
$(h_1([X,Y;\theta]),h_{2,E}((Y_E,\xi\circ\theta_E^{-1})))$. \hqed
\end{remark}

\begin{definition} \label{locfree} A finitely generated $\A$-module
$M$ is said to be \textit{locally free} if $M_v$ is a free
$\A_v$-module for all places $v$ in $S_f(F)$. It follows easily from
this definition that if $M$ is a locally free $\A$-module, then it is
projective and the associated $A$-module $M_F$ is free. \hqed
\end{definition}

\begin{remark}\label{exdes} If $[X,Y;\xi]$
is an element of $K_0(\A,E)$ for which both $X$ and $Y$ are locally
free $\A$-modules, then one can give an explicit representative for
$h_{\A,E}([X,Y;\xi])$ as follows.

For any ordered set of $d$ elements $\{e^j: 1 \le j \le d\}$ we
write $\underline{e^j}$ for the $d \times 1$ column vector with
$j$-th entry equal to $e^j$.

We choose an $A$-basis $\{ y^j\}$ of $Y_F$ and, for each $v\in
 S_f(F)$, an $\A_v$-basis $\{y_v^j\}$ of $Y_v$ and an $\A_v$-basis
 $\{x_v^j\}$ of $X_v$. Let $\mu_v$ be the element of $GL_d(A_v)$ such
 that $\underline{y_v^j} = \mu_v \underline{y^j}$. We choose $\theta \in
 \Is_A(X_F,Y_F)$. Since $\{\theta^{-1}(y^j)\}$ is an $A$-basis of $X_F$, there
 exists an element $\lambda_v$ of $GL_d(A_v)$ such that
 $\underline{x_v^j} = \lambda_v\underline{ \theta^{-1} (y^j)}$. Finally we
 let $\mu \in GL_d(A_E)$ denote the matrix of $\xi\circ(\theta^{-1}\otimes
 _FE)$ with respect to the $A$-basis $\{ y^j\}$ of $Y_E$.  Then $h_{\A
 ,E}([X,Y;\xi])$ is represented by the function
\[ 
\left( \prod_{v \in S_f(F)} \Det(\lambda_v\mu_v^{-1}) \right) \times
 \Det (\mu)\in \Hom (R_\Gamma ,J_f(\bQ^c))^{\Omega_F}\times
 \Hom(R_\Gamma , (E^c)^{\times}).
\]
This construction will be used in the proof of Proposition \ref{pr}. \hqed
\end{remark}

\begin{examples}\label{exdes1} Suppose that $F$ is a number field, and
let $R=O_F$.

(i) Suppose that $\A = R[G]$ (see Example \ref{class}). Then condition
(\ref{hd}) is satisfied if we take $\Gamma = G$, viewed as a trivial
$\Omega_F$-module.

(ii) In Example \ref{geom}, set $\Gamma = G(F^c)$, endowed with the
natural $\Omega_F$-action. Then $G \otimes_R F^c$ is a finite constant
group scheme over $\Spec(F^c)$, and so its Cartier dual $G^* \otimes_R
F^c$ is equal to $\Spec(F^c[\Gamma])$ (see e.g. \cite[Section
2.4]{W1}). Hence $\A = \calO_{G^*}$ satisfies condition \eqref{hd}.
Furthermore, in this case the description of Theorem \ref{Hom} may be
interpreted more explicitly as follows.

Let $J_f(A)$ denote the group of finite ideles of $A$, i.e. if $\M$ is
the (unique) maximal $R$-order in $A$, then $J_f(A)$ is the restricted
direct product of the groups $A_{v}^{\times}$ with respect to the
subgroups $\M_{v}^{\times}$ for all places $v \in S_f(F)$. Define
a map
\begin{equation*}
\Delta '_{\A ,E}: A^{\times} \to \frac{J_f(A)}{U_f(\A)} \times
A_E^{\times}; \qquad a \mapsto (a) \times a^{-1}.
\end{equation*}
Then, taken in conjunction with Remark \ref{ct}, the result of
Theorem \ref{Hom} implies that there is a natural isomorphism
\begin{equation}\label{E:idelic}
K_0(\A, E) \xrightarrow{\sim} \Cok(\Delta '_{\A ,E}).
\end{equation}
\hfill \hqed
\end{examples}

\section{Metrised structures} \label{S:metrised}

In this section we let $R$ denote the ring of algebraic integers
in a number field $F$. We fix an $R$-order $\A$ in a finite
dimensional $F$-algebra $A$, and we continue to assume that the
condition (\ref{hd}) is satisfied.

For each $\phi \in \hat\Gamma$ we write $W_\phi$ for the Wedderburn
component of $F^c[\Gamma]$ which corresponds to the contragredient
character $\overline{\phi}$ of $\phi$. For any $F^c[\Gamma ]$-module
$X$ we then set
\begin{equation*}
X_\phi := \bigwedge^{\top}_{F^c}((X\otimes_{F^c}W_\phi)^\Gamma),
\end{equation*}
where `$\bigwedge^{\top}_{F^c}$' denotes the highest exterior power
over $F^c$ which is non-zero, and $\Gamma$ acts diagonally on
$X\otimes_{F^c}W_\phi$.

For each $v \in S_\infty (F)$ we fix an embedding $\sigma_v :F^c
\rightarrow F_v^c$ (which induces the place $v$ upon restriction to
$F$), and we also identify $F_v^c$ with $\bC$.

For each complex number $z$ we write $\overline{z}$ for its complex
conjugate.

\subsection{Metrised $\A$-modules}
\medskip

\begin{definition} Let $X$ be a
finitely generated locally-free $\A$-module. A {\it metric} $\mu$ on
$X$ is a set $\{ \mu_{v,\phi}: v \in S_\infty(F), \phi \in
\hat\Gamma\}$ where, for each $v \in S_\infty(F)$ and $\phi \in
\hat\Gamma$, $\mu_{v,\phi}$ is a hermitian metric on the complex line
$(X\otimes_{R}F^c)_\phi \otimes_{F^c,\sigma_v}\bC$. If $x \in
(X\otimes_{R}F^c)_\phi \otimes_{F^c,\sigma_v}\bC$, then we shall write
$\mu_{v,\phi}(x)$ for the length of $x$ with respect to the metric
$\mu_{v,\phi}$.
 
A {\it metrised $\A$-module} is a pair $(X,\mu)$ consisting of a
locally-free $\A$-module $X$ and a metric $\mu$ on $X$. \hqed
\end{definition}

\begin{example}\label{exex} {\rm (Classical) In this example we adopt
the notation of Example \ref{class}. We take $\A$ to be $\bz [G]$, so
that $\Gamma = G$ (viewed as a trivial $\Omega_{\bQ}$-module). Via the
fixed embedding $\sigma_\infty:\bQ^c \to \bC$ restricted to the normal
closure $N$ of $L$, we identify $\Sigma_\bq(L,N)$ with the set
$\Sigma(L)$ of embeddings of $L$ into $\bc$. We remind the reader that
if $X$ is a finitely generated, locally free $\fA$-module, then any
$G$-equivariant positive definite hermitian form on $X_{\bC}$ induces
a metric on $X$ in a natural way (see e.g. \cite[Definition
2.2]{CPT}).

We write $\mu_{\bC [G]}$ for the $G$-equivariant positive definite
hermitian form on $\bC [G]$ which satisfies
\[ 
\mu_{\bC [G]}( \sum_{g \in G}x_gg, \sum_{h \in G}y_hh) =
\sum_{g\in G}x_g\overline{y_g}.
\]

(i) There is a $G$-equivariant positive definite hermitian form
$\mu_L$ on $H^L_\bQ \otimes_\bZ \bC$ which is defined by the rule
\[ \mu_L(
\sum_{\sigma \in \Sigma(L)}z_\sigma ', \sum_{\tau \in
\Sigma(L)}z_\tau) = \sum_{\sigma \in \Sigma (L)}z_\sigma
'\overline{z_\sigma}.
\] 
For each $\phi \in \hat G$, we write $\mu_{L,\infty,\phi}$ for the
metric on $(H^{L}_{\bQ} \otimes_{\bZ} \bQ^c)_\phi$ that is obtained as
the highest exterior power of the metric on
$$
(H^{L}_{\bQ} \otimes_{\bZ} W_{\phi})^G =
(H^{L}_{\bQ} \otimes_{\bZ} \bC) \otimes_{\bC} (W_{\phi}
\otimes_{\bQ^c, \sigma_\infty} \bC)^G
$$
which is induced by $\mu_L$ on $H^{L}_{\bQ} \otimes_{\bZ} \bC$ and by
the restriction of $\mu_{\bc [G]}$ on $W_{\phi} \otimes_{\bQ^c,
\sigma_\infty} \bC$.
\smallskip

(ii) There is a $G$-equivariant positive definite hermitian form $h_L$
on $\Ce\otimes_\bq L$ which is defined by the rule
\[ h_L(z_1\otimes m,z_2\otimes n) = z_1\overline{z_2}\sum_{\sigma \in
\Sigma(L)}\sigma(m) 
\overline{\sigma (n)} .\] We recall that in \cite[\S5.2]{CPT} this
form is referred to as the `Hecke form'.

For each $\phi \in \hat G$ we write $h_{L,\infty,\phi}$ for the metric
on $(L \otimes_{\bQ} \bQ^c)_\phi$ that is obtained as the highest
exterior power of the metric on
$$
(L \otimes_{\bQ} W_{\phi})^G \otimes_{\bQ^c, \sigma_\infty} \bC =
((\bC \otimes_{\bQ} L) \otimes_{\bC} (W_{\phi} \otimes_{\bQ^c,
\sigma_\infty} \bC))^G
$$
which is induced by $h_L$ on $\bC \otimes_{\bQ} L$ and by the
restriction of $\mu_{\bc [G]}$ on $W_{\phi} \otimes_{\bQ^c,
\sigma_\infty} \bC$.

We set $h_{L,\bullet} := \{ h_{L,\infty,\phi }:\phi \in \hat G\}$. If
$\calL$ is any full projective $\bz [G]$-sublattice of $L$, then
the pair $(\calL, h_{L, \bullet})$ is a metrised $\bz [G]$-module.} \hqed
\end{example}

\begin{example} \label{exex2} (Geometrical) In this example we
adopt the notation of Example \ref{geom}.
\smallskip

(i) For each $v \in S_\infty(F)$ and $\phi\in \hat\Gamma$, the space
$(A\otimes_FF^c)_\phi\otimes_{F^c,\sigma_v}\bc$ identifies naturally
with $\bc$ and so is endowed with a metric $\mu_{A,v,\phi}$ coming
from the standard metric on $\bc$. The set $\mu_A := \{
\mu_{A,v,\phi}: v \in S_\infty (F),\phi \in \hat\Gamma\}$ is then a
metric (`the trivial metric') on $\A$.
\smallskip

(ii) It is shown in \cite[\S2]{AP} that, for each place $v \in
S_\infty (F)$ and character $\phi \in \hat\Gamma$, the space
$(\calL_\pi\otimes_RF^c)_\phi\otimes_{F^c,\sigma_v}\bc$ may be endowed
with a canonical `Neron metric' $||\cdot ||_{v,\phi}$ which is
constructed by using canonical splittings of the extension
\eqref{E:extension}. (These canonical splittings are analogous to the
canonical splittings of extensions of abelian varieties by tori that
may be used to define N\'eron pairings on abelian varieties (see
e.g. \cite{MT}).)  Let $||\cdot ||_\pi$ denote the family of metrics
$\{ ||\cdot ||_{v,\phi}: v\in S_\infty (F), \phi \in
\hat\Gamma\}$. Then the pair $(\calL_\pi,||\cdot ||_\pi)$ constitutes
a metrised $\A$-module.  \hqed
\end{example}

\subsection{Arithmetic class groups} In this subsection we use idelic
valued functions on $R_\Gamma$ to define a group which classifies the
structure of metrised $\A$-modules. All of the constructions in this
subsection are motivated by those of \cite[\S3.1, \S3.2]{CPT}.
\smallskip

We write $|\Delta |_\A$ for the homomorphism
\begin{equation*}\label{d2} 
\Det (A^\times) \to \frac{\Hom
(R_\Gamma,J_f(F^c))^{\Omega_F}}{\Det(U_f(\A ))} \times \Hom
\left(R_\Gamma,\prod_{v \in S_\infty (F)}\br^\times_{>0}\right);
\qquad \theta \mapsto (\theta, |\theta|)
\end{equation*}
 where $|\theta|$ is the homomorphism which sends each $\phi\in
\hat\Gamma$ to the element
\[
\prod_{v \in S_\infty(F)}|\sigma_v(\theta(\phi))|^{-1} \in \prod_{v \in
 S_\infty(F)}\br^\times_{>0}.
\]

\begin{remark} \label{R:caveat} 
We caution the reader that the definition of the homomorphism
$|\theta|$ given above is slightly different from that used in
\cite[\S3.1]{CPT}. Our $|\theta|(\phi)$ is equal to
$|\theta|(\phi)^{-1}$ in the notation of \cite{CPT}. This reflects our
definitions of $\Delta_{\A,E}$ in \eqref{E:deltadef} and $h_{\A,E}$ in
\eqref{E:homiso}. Our normalisation is also consistent with the
definition of the Hermitian classgroup given in \cite[Chapter II,
\S5]{Fr}.  \hqed
\end{remark}

\begin{definition} \label{D:AC}
We define the {\it arithmetic classgroup} $\AC(\A)$
of $\A$ to be the cokernel of $|\Delta|_\A$. \hqed
\end{definition}

\begin{remark}
When $F=\Q$ and $\A = \Z[\Gamma]$ for some finite group $\Gamma$, then
it is easy to show (taking into account Remark \ref{R:caveat} above)
that $\AC(\A)$ is isomorphic to the arithmetic classgroup $A(\A)$
which is described in \cite[Definition 3.2]{CPT}. \hqed
\end{remark}

We now suppose given a metrised $\A$-module $(M ,\mu)$.  Let $M$ have
rank $d$ over $\A$, and choose an $A$-basis $\{ m^j\}$ of $M_F$ and,
for each $v \in S_f(F)$, an $\A_v$-basis $\{m_v^j\}$ of $M_v$.  Then
there exists an element $\lambda_v$ of $GL_d(A_v)$ such that
$\underline{m_v^j} = \lambda_v\underline{m^j}$.

For each element $m$ of $M_F$ we set
\[ 
r(m) := \sum_{\gamma \in \Gamma}\gamma m\otimes \gamma \in
M \otimes _{R} F^c[\Gamma].
\] 
For each element $w$ of $W_\phi$ we have
\[ 
r(m)(1\otimes w) \in (M \otimes_{R}W_\phi)^\Gamma.
\]
Let $\{ w_{\phi,k}:k\}$ be an $F^c$-basis of $W_\phi$ which is
orthonormal with respect to the restriction of $\mu_{F^c [\Gamma]}$ to
$W_\phi$. Then the set $\{ r(m^j)(1\otimes w_{\phi,k}): j, k\}$ is an
$F^c$-basis of $(M \otimes_{R} W_\phi )^\Gamma$, and so the wedge
product
\[ 
\bigwedge_j \bigwedge_k r(m^j)(1\otimes w_{\phi,k})
\]
is an $F^c$-basis of the line $(M \otimes_{R}F^c)_\phi$.

\begin{definition} We define $[M, \mu]$ to be the element of
$\mathrm{AC}(\A)$ which is represented by the homomorphism on
 $R_\Gamma$ which sends each character $\phi \in \hat \Gamma$ to
\[ 
\prod_{v\in S_f(F)}\Det(\lambda_v)(\phi)\times \prod_{v \in
 S_\infty (F)}\mu_{v,\phi} 
\left( \left(\bigwedge_j \bigwedge_k r(m^j)(1\otimes
w_{\phi,k})\right)\otimes_{F^c,\sigma_v}1\right)^{\frac{1}{\phi (1)}}\in
J_f(F^c)\times \prod_{v \in S_\infty (F)}\br^\times_{>0} .
\] 
It can be shown that $[M,\mu]$ is independent of the precise choices
of bases $\{ m^j\}$, $\{m_{v}^{j}\}$ and $\{w_{\phi,k}\}$. \hqed
\end{definition}

\subsection{The connection to generalised twisted forms} \label{S:concat} 
Let $E$ be any field which contains $F^c$ and which is such that, for
each $v \in S_\infty (F)$, the fixed embedding $\sigma_v:F^c \to \bC$
factors through an embedding $\tilde \sigma_v:E \to \bC$.

The map
\[ 
\pi_E: E^\times \to \prod_{v
\in S_\infty (F)}\br^\times_{>0}; \qquad e \mapsto \prod_{v \in
S_\infty (F)} |\tilde\sigma_v(e)|
\] 
induces a homomorphism
\begin{equation*}
\varepsilon: \Det(A_{E}^{\times}) \simeq \Hom(R_{\Gamma},E^{\times})
\longrightarrow
\Hom(R_{\Gamma}, \prod_{v \in S_\infty (F)}\br^\times_{>0}).
\end{equation*}
Then, via the Hom-descriptions of Theorem \ref{Hom} and Definition
\ref{D:AC}, it is not hard to see that the map
\begin{equation*}
\delta: \frac{\Hom(R_\Gamma ,J_f(F^c))^{\Omega_F}}{\Det(U_f(\A ))}
\times \Det(A_{E}^{\times}) \longrightarrow
\frac{\Hom(R_\Gamma ,J_f(F^c))^{\Omega_F}}{\Det(U_f(\A ))}
\times 
\Hom(R_{\Gamma}, \prod_{v \in S_\infty (F)}\br^\times_{>0})
\end{equation*}
given by $x_1 \times x_2 \mapsto x_1 \times \varepsilon(x_2)$ induces
a homomorphism
\begin{equation*}
\partial_{\A,E}:K_0(\A,E) \to \AC(\A).
\end{equation*}

\begin{lemma} There is a natural isomorphism
\begin{equation} \label{E:circle}
\ker(\partial_{\A,E}) \xrightarrow{\sim} \frac{\Hom(R_\Gamma
,\ker(\pi_E))}{[\Hom(R_\Gamma ,\ker(\pi_E)) \cap \Det(A^{\times})]
\cap \Det( U_f(\A))}. 
\end{equation}
(In the denominator of the right-hand side of \eqref{E:circle}, the
first intersection takes place in 
$$
\Hom(R_{\Gamma}, E^{\times}) \simeq
\Det(A_{E}^{\times}),
$$
and then the second intersection takes place in
$\Hom(R_{\Gamma},J_f(F^c))$.)
\end{lemma}
\begin{proof}
Set 
\begin{equation*}
N_1 :=  \frac{\Hom(R_\Gamma ,J_f(F^c))^{\Omega_F}}{\Det(U_f(\A ))}
\times \Det(A_{E}^{\times}),
\end{equation*}
and write $N_2:= \delta(N_1)$. Then there is a commutative diagram:
\begin{equation*}\minCDarrowwidth1em
\begin{CD} 
0@ > >> \Hom(R_\Gamma,\ker(\pi_E)) @>>> N_1 @>
\delta >> N_2 @> >> 0\\ 
@. @. @A\Delta_{\A ,E}AA @A|\Delta |_\A AA \\ @. @.
 \Det (A^\times) @= \Det( A^\times).
\end{CD}
\end{equation*}
Applying the Snake Lemma to this diagram (and taking into account
Theorem \ref{Hom} and Definition \ref{D:AC}) yields the sequence
$$
\Hom(R_\Gamma, \Ker(\pi_E)) \to K_0(\fA,E)
\xrightarrow{\partial_{\fA,E}} \mathrm{AC}(\fA) \to 0,
$$
and it is easy to check (again using Theorem \ref{Hom}) that this
implies the desired result.
\end{proof}
\smallskip

Given an explicit element of $K_0(\A, E)$ we now aim to describe
its image under $\partial_{\A,E}$.
\smallskip

\begin{definition} Let $X_1$ and $X_2$ be finitely
generated locally-free $\A$-modules and suppose that $\xi \in
\Is_{A_E}(X_{2,E}, X_{1,E})$. For each place $v \in
S_\infty (F)$ and character $\phi \in \hat\Gamma$ we write 
\[
\xi_{v,\phi} : (X_2\otimes_{R}F^c)_\phi\otimes_{F^c,\sigma_v}\bc
\xrightarrow{\sim}
(X_1\otimes_{R}F^c)_\phi\otimes_{F^c,\sigma_v}\bc
\] 
for the isomorphism of complex lines which is induced by $\xi$. If $h$
is any metric on $X_1$, then the \textit{pullback} of $h$ under $\xi$ is
defined to be the metric $\xi^*(h)$ on $X_2$ which satisfies
\[ 
\xi^*(h)_{v,\phi}(z) = h_{v,\phi}(\xi_{v, \phi}(z))
\]
for each $v \in S_\infty (F), \phi\in \hat\Gamma$ and $z \in
(X_2\otimes_{R} F^c)_\phi\otimes_{F^c,\sigma_v}\bc$. \hqed
\end{definition}

\begin{example}\label{exex3}{\rm (i) (Classical) An explicit
computation shows that, in the notation of Example \ref{exex}, one has
$(\pi^{L,N}_\bq)^{*}(\mu_{L,\bullet}) = h_{L,\bullet}$.
\smallskip

(ii) (Geometrical) We use the notation of Example \ref{exex2}. If
$\xi_\pi$ is any isomorphism as in (\ref{E:splitting}), then it
follows from \cite[Remark 2.3]{AP} that $\xi_{\pi}$ induces an
isometry between $\L_{\pi}$ endowed with the N\'eron metric $||\cdot
||_\pi$ and $\A$ endowed with the trivial metric. (Note that the map
which is denoted by $\xi_{\pi}$ in \cite{AP} is equal to
$\xi_{\pi}^{\otimes e}$ in our present notation, where $e$ denotes the
exponent of the group scheme $G$.) Hence we have that
$(\xi_\pi\otimes_{R^c}F^c)^*(\mu_A)=||\cdot ||_\pi.$} \hqed
\end{example} 

\begin{proposition}\label{pr} Let $(X,Y;\xi)$ be any object of
$\calP_\A\times_E \calP_{\A}$ for which $X$ and $Y$ are both locally
free $\A$-modules. Then, for any metric $\rho$ on $Y$, one has
\[ 
\partial_{\A,E}(h_{\A,E}([X,Y;\xi])) = [X,\xi^*(\rho)] - [Y,\rho]
\in \mathrm{AC}(\A).
\] 
\end{proposition}

\begin{proof} We use the notation introduced in Remark \ref{exdes}.

It is clear that $\partial_{\A,E} (h_{\A,E}([X,Y;\xi]))$ is
represented by the homomorphism which sends each character $\phi
\in \hat \Gamma$ to
\begin{equation}\label{fo}
\prod_{v\in S_f(F)}\Det(\lambda_v\mu_v^{-1})(\phi) \times
\prod_{v\in S_\infty(F)}|\tilde\sigma_v(\Det
(\mu)(\phi))|.
\end{equation}

For each element $x$ of $X_{i,\phi} := (X_i\otimes_{R}F^c)_\phi
\otimes_{F^c}E$ (where $i=1,2$) and each place $v \in S_{\infty}(F)$, we
write $[x]_v$ for the image of $x$ in $X_{i,\phi}\otimes_{E,\tilde
\sigma_v}\bc =
(X_i\otimes_{R}F^c)_\phi\otimes_{F^c,\sigma_v}\bc$.

Using this notation, the element $[X,\xi^*(\rho)]$, resp.
$[Y,\rho]$, is represented by the homomorphism which sends each
character $\phi \in \hat \Gamma$ to
\begin{equation}\label{ft}
\prod_{v\in S_f(F)}\Det (\lambda_v)(\phi) \times \prod_{v \in
S_\infty(F)}\rho_{v,\phi}\left( \left[\bigwedge_j \bigwedge_k r(
\xi\circ \psi (y^j))(1\otimes w_{\phi
,k})\right]_v\right)^{\frac{1}{\phi(1)}},
\end{equation}
resp.
\begin{equation}\label{ftt}
\prod_{v\in S_f(F)}\Det (\mu_v)(\phi) \times \prod_{v \in
S_\infty(F)}\rho_{v,\phi}\left( \left[\bigwedge_j \bigwedge_k
r(y^j)(1\otimes w_{\phi ,k})\right]_v\right)^{\frac{1}{\phi(1)}}.
\end{equation}
For each $z = \sum_{\gamma \in \Gamma}c_\gamma \gamma \in E[\Gamma]$
we set $\overline{z} := \sum_{\gamma \in \Gamma}c_\gamma \gamma^{-1}$,
and we extend this convention to matrices over $E[\Gamma ]$ by
applying it to individual entries. Now
\[
\xi \circ \psi (y^j) = \sum_{l}\mu_{l,j}y^l
\]
where $\mu_{l,j}$ is the $(l,j)$-component of the matrix $\mu$,
and so
\[ 
r(\xi\circ \psi (y^j))(1\otimes w_{\phi, k}) = \sum_l
r(y^j)(1\otimes \overline{\mu_{l,j}}w_{\phi,k}).
\] 
This implies that
\[ 
\bigwedge_j \bigwedge_k r(\xi\circ \psi (y^j))(1\otimes w_{\phi,k}) =
 \Det(\overline{\mu })(\overline{\phi})^{\phi(1)}\cdot
   \bigwedge_j \bigwedge_k r(y^j)(1\otimes w_{\phi,k}).
\]
But $\Det(\overline{\mu })(\overline{\phi}) = \Det(\mu (\phi))$,
and so for each $v \in S_\infty(F)$ one has
\[
\frac{\rho_{v,\phi}\left( \left[\bigwedge_j \bigwedge_k  r( \xi \circ \psi
(y^j))(1\otimes w_{\phi
,k})\right]_v \right)^{\frac{1}{\phi(1)}}}{\rho_{v,\phi}\left(
\left[\bigwedge_j \bigwedge_k r( 
y^j)(1\otimes w_{\phi ,k})\right]_v\right)^{\frac{1}{\phi(1)}}} = |\tilde
\sigma_v(\Det(\mu) (\phi))| \in \br^\times_{> 0}.
\] 
It is now clear that the expression (\ref{ft}) is equal to the product
of the expressions (\ref{fo}) and (\ref{ftt}), and this immediately
implies the claimed equality.
\end{proof}

\section{The classical case} \label{S:classical}

In this section we give two applications of our approach in the
context of Example \ref{class}. We thus fix a finite tamely ramified
Galois extensions of number fields $L/K$, with $F \subseteq K$, and we
set $G := \Gal(L/K)$. Our interest is in the element
\[ 
\delta^L_{O_F[G]}(O_L) := [O_L,H_F^L\otimes_\bz O_F;\pi^{L,N}_F] \in
K_0(O_F [G],N).
\]

\subsection{Arithmetic classes} In this subsection we combine
 Proposition \ref{pr} with a result of Bley and the second named
author in \cite{BB} in order to describe an explicit homomorphism
which represents the element $\delta^{L}_{\bZ[G]}(O_L)$. We then explain
how this yields an explicit homomorphism which represents the element
$[O_L,h_{L,\bullet}]$ of $\mathrm{AC}(\bZ [G])$. We show that this
description of $[O_L,h_{L,\bullet}]$ gives a natural refinement of a
result of Chinburg, Pappas and Taylor in \cite{CPT}.

Before stating our explicit description of the class
$\delta^{L}_{\bZ[G]}(O_L)$ we must introduce some notation.
 
For each character $\phi \in \hat G$, we write $\tau(K,\phi)$ for the
Galois-Gauss sum associated to $\phi$, and we let $\epsilon(K ,\phi)$
denote the epsilon constant which arises in the functional equation of
the Artin $L$-function attached to $\phi$ (see \cite[Chapter I,
\S5]{Fr0}). With respect to the canonical identification $\zeta (\bc
[G]) = \prod_{\hat G}\bc$ (which is induced by the fixed embedding
$\sigma_\infty: \bQ^c \to \bC$) we then set
\[ 
\epsilon_{L/K} :=
(\epsilon (K,\overline{\phi }))_\phi \in \zeta (\bc [G])^\times.
\]

This element is the natural epsilon constant which is associated to
the $\zeta ( \bc [G])$-valued equivariant Artin $L$-function of $L/K$,
and it actually lies in $\zeta (\br [G])^\times$. The last fact
implies that we may therefore choose an element
\begin{equation*}
\lambda = (\lambda_\phi)_\phi \in \zeta(\bQ[G])^{\times} \subset
\zeta(\bC[G])^{\times}
\end{equation*}
such that
\begin{equation} \label{pc}
\lambda \epsilon_{L/K} \in \Im (\nr_{\br [G]})
\end{equation}
(see \cite[\S 3.1]{BB}). We write $d_K$ for the absolute discriminant
of $K$.
 
We write $\iota: N \rightarrow \bC$ for the embedding induced by
restricting $\sigma_\infty$ to $N$, and we use this to view
$\delta^{L}_{\bZ[G]}(O_L)$ as an element of $K_0(\bZ[G],\bC)$ (using
the inclusion $K_0(\bZ[G],N) \subseteq K_0(\bZ[G],\bC)$ induced by the
functor $-\otimes_{N,\iota}\bc$).

We can now state the main result of this section.

\begin{theorem} \label{kcr}  Let $\lambda$ be any element of
$\zeta(\bq [G])^\times$ which satisfies \eqref{pc}. Then the class
$\delta^{L}_{\bZ[G]}(O_L)$ in $K_0(\bZ[G],\bC)$ is represented by the
homomorphism which sends each character $\phi \in \hat G$ to
\begin{equation}\label{fs}
\lambda_\phi^{-1}\prod_{v \in S_f(K)}y_v(\phi) \times
\lambda_\phi\tau(K,\phi) |d_K|^{\frac{\phi (1)}{2}}.
\end{equation}
Here $y_v$ denotes the `unramified characteristic' function
on $R_G$ defined in \cite[Chapter IV, \S1]{Fr0}.
\end{theorem}

We shall derive a description of the class $\delta_{\bZ[G]}^{L}(O_L)$
from certain results in \cite{BB}. For the reader's convenience we
first recall some notation from loc. cit.

\begin{definition} \label{eqdisc} (See \cite[\S3.2]{BB}) Set $H_L:=
\prod_{\Sigma (L)} \bZ$. The groups $G$ and $\Gal(\bC/\bR)$ act on
$\Sigma(L)$, and they endow $H_L$ with the structure of $G \times
\Gal(\bC/\bR)$-module. For any $\Gal(\bC/\bR)$-module $X$, write $X^+$
and $X^-$ for the submodules on which complex conjugation acts by $+1$
and $-1$ respectively.
 
There is a canonical $\bC[G \times \Gal(\bC/\bR)]$-equivariant
isomorphism
\begin{equation*}
\rho_L: L \otimes_{\bQ} \bC \xrightarrow{\sim} H_L \otimes_{\bZ} \bC;
\qquad l \otimes z \mapsto (\sigma(l)z)_{\sigma \in \Sigma(L)}.
\end{equation*}
 
We write $\pi_L$ for the composite $\bR[G]$-equivariant isomorphism
given by
\begin{align*}
L \otimes_{\bQ} \bR & = (L \otimes_{\bQ} \bC)^+ \\ &
\xrightarrow{\rho_L} (H_L \otimes_{\bZ} \bC)^+ \\ &= (H_{L}^{+}
\otimes_{\bZ} \bR) \oplus (H_{L}^{-} \otimes_{\bZ} (i \cdot \bR))
\\ &\xrightarrow{\times (1,-i)} (H_{L}^{+} \otimes_{\bZ} \bR) \oplus
(H_{L}^{-} \otimes_{\bZ} \bR) \\ &= H_L \otimes_{\bZ} \bR.
\end{align*}

We now define the \textit{$G$-equivariant discriminant}
$\delta_{L/K}(O_L)$ of $O_L$ by setting
\begin{equation*}
\delta_{L/K}(O_L):= [O_L,H_L; \pi_L] \in K_0(\bZ[G],\bR).
\end{equation*}
In what follows, we shall in fact view $\delta_{L/K}(O_L)$ as an
element of $K_0(\bZ[G],\bC)$ via the natural inclusion
$K_0(\bZ[G],\bR) \subseteq K_0(\bZ[G],\bC)$.\hqed
\end{definition}
\bigskip

\textit{Proof of Theorem \ref{kcr}.}
It follows from results in \cite[see especially Lemma 7.4 and
Corollary 7.7]{BB} that the class $\delta_{L/K}(O_L)$ is represented
in the Hom-description of Theorem \ref{Hom} by the element of
$\Hom(R_G,J_f(\bq^c))^{\Omega_\bq}\times \Hom(R_G,\bc^\times)$ which
sends each character $\phi \in \hat G$ to
\[
\lambda_\phi^{-1}\prod_{v \in
  S_f(K)}y_v(\phi) \times\lambda_\phi\epsilon (K,\phi),
\]
where $y_v$ is as defined in the statement of Theorem \ref{kcr}.
 
For each place $v \in S_{\infty}(K)$, we let $G_v$ denote the
decomposition group of $v$ in $G$. If $\phi \in \hat G$, we let
$V_{\phi}$ be a representation space for $\phi$, and we define a
complex number $w_{\infty}(K,\phi)$ by setting
\begin{equation*}
w_{\infty}(K,\phi):= \prod_{v \in S_{\infty}(K)}
i^{-\codim(V_{\phi}^{G_v})}.
\end{equation*}
 
We may identify $H_L$ with $H_\bQ^L$ via our fixed embedding $\iota:
N \rightarrow \bC$. By comparing the homomorphisms
$\pi_{\bQ}^{L,N}\otimes_{N, \iota}\bC$ and $\pi_L\otimes_\bR\bC$, it
may be shown that the element $\delta_{L/K}(O_L) -
\delta^{L}_{\bZ[G]}(O_L)$ of $K_0(\bZ[G],\bC)$ is represented by the
element of $\Hom(R_G,J_f(\bq^c))^{\Omega_\bq} \times
\Hom(R_G,\bc^\times )$ which sends each character $\phi \in \hat G$ to
$1 \times w_\infty(K,\overline{\phi})$.
 
In addition, for each $\phi \in \hat G$ one has
\begin{equation}\label{et}
\epsilon(K,\phi) =
\tau(K,\phi)w_\infty(K,\overline{\phi})|d_K|^{\frac{\phi (1)}{2}}
\end{equation}
(see \cite[(13) and the displayed formula which follows it]{BB}).
 
It now follows that $\delta_{\bz
[G]}^L(O_L)$ is represented by the homomorphism which sends each
character $\phi \in \hat G$ to
\begin{equation*}
\lambda_\phi^{-1}\prod_{v \in S_f(K)}y_v(\phi) \times
\lambda_\phi\tau(K,\phi) |d_K|^{\frac{\phi (1)}{2}},
\end{equation*}
as claimed. \hqed

\begin{corollary} \label{C:kcr} 
Let $\lambda$ be any element of $\zeta(\bq [G])^\times$ which
satisfies \eqref{pc}. Then the class $[O_L, h_{L,\bullet}]$ in
$\AC(\bZ[G])$ is represented by the homomorphism which sends each
character $\phi \in \hat G$ to 
\[
\lambda_\phi^{-1} \times
|\lambda_\phi\tau(K,\phi)|(|G|^{[K:\bq]}|d_K|)^{\frac{\phi(1)}{2}}.
\]
\end{corollary}
 
\begin{proof} We first observe that Example \ref{exex3}(i) and
Proposition \ref{pr} imply that 
\begin{equation} \label{ss}
[O_L,h_{L,\bullet}] = [H^L_\bq,\mu_{L,\bullet}] + \partial_{\bz
[G],\bc}(\delta_{\bz [G]}^L(O_L))   \in \mathrm{AC}(\bz [G]).
\end{equation}

It follows via an explicit computation (using, for example,
\cite[Lemma 2.3]{CPT}) that the class $[H^L_\bq,\mu_{L,\bullet}]$ is
represented by the homomorphism which sends each character $\phi\in
\hat G$ to $1 \times |G|^{[K:\bq]{\frac{\phi(1)}{2}}}$. Combining this
with \eqref{ss} and Theorem \ref{kcr}, we deduce that 
$[O_L,h_{L,\bullet}]$ is represented by the homomorphism which sends
 each $\phi \in \hat G$ to
\[ 
\lambda^{-1}_\phi\prod_{v \in S_f(K)}y_v(\phi) \times
|\lambda_\phi\tau(K,\phi)|(|G|^{[K:\bq ]} |d_K|)^{\frac{\phi (1)}{2}}.
\] 
The desired result now follows from this description upon
noting that, as a consequence of \cite[Chapter IV, Theorem 29(i)]{Fr0},
the homomorphism which sends each $\phi \in \hat G$ to
$$
\left(\prod_{v \in S_f(K)}y_v(\phi)\right)\times 1
$$ 
belongs to the image of $|\Delta |_\A$.
\end{proof}

It is often convenient to replace $\AC(\bZ[G])$ by weaker classifying
groups. As an example of this, we shall now explain how Corollary
\ref{C:kcr} yields a natural refinement of \cite[Theorem 5.9]{CPT}.

We first quickly recall the definition of the `tame symplectic
arithmetic classgroup' $A^s_T(\bz [G])$ of $\bZ[G]$ from [loc. cit.,
\S4.3]. We write $R_G^s$ for the subgroup of $R_G$ which is generated
by the irreducible symplectic $F^c$-valued characters of $G$, and we
write $\Det ^s(-)$ for the restriction of the function $\Det(-)$ to
$R_G^s$. Let $T$ denote the maximal abelian tamely ramified extension
of $\bq$ in $\bq^c$, and set $\hat \bZ= \varprojlim_n \bZ/n\bZ$. Then
we write $\Det (\hat O_T[G]^\times)$ for the direct limit of $\Det
(\hat O_N [G]^\times)$, where $N$ runs over all finite extensions of
$\bq$ in $T$ and $\hat O_N$ denotes the ring of integral adeles
$\hat\bz \otimes O_N$.

\begin{definition} \label{D:sympclass}
The tame arithmetic symplectic classgroup $A^s_T(\bz [G])$ of $\bZ[G]$
is defined to be the cokernel of the
map
\[
\Det^s(\bq [G]^\times) \lra \frac{\Det^s(\hat
O_T[G]^\times)\Hom(R_G^s,J_f(\bq^c))^{\Omega_\bq}} {\Det^s(\hat
O_T[G]^\times)}\times \Hom(R_G^s,\br_{>0}^\times )
\]
which is induced by the diagonal map
\[
\Delta^s : \Det^s(\bq [G]^\times) \lra
\Hom(R_G^s,J_f(\bq^c))^{\Omega_\bq}\times \Hom(R_G^s,\br_{>0}^\times
);\,\,\,\, \theta \mapsto (\theta, \theta^{-1}).
\]
\hqed
\end{definition}

We write
\[
\theta^s_T: \Hom(R_G^s,J_f(\bq^c))^{\Omega_\bq}\times
\Hom(R_G^s,\br_{>0}^\times ) \lra A^s_T(\bz [G])
\]
for the natural surjective map, and
\[
\pi^s_T: \mathrm{AC}(\bz [G]) \lra A^s_T(\bz [G])
\]
for the homomorphism induced by restricting an element of
$\Hom(R_G,J_f(\bq^c))^{\Omega_\bq}\times\Hom(R_G,\br_{>0}^\times)$ to
$R_G^s$, and then applying $\theta_T^s$.

\begin{definition} \label{D:Pfaffian}
Let $\mathrm{Pf}(O_L)$ denote the `Pfaffian' element of $\Hom
(R_G^s,J_f(\bq^c ))$ which is defined (in [loc. cit., Def. 5.8]) as
follows: for each $\phi \in R_G^s$ one has
\[
\mathrm{Pf}(O_L)(\phi)_p = \prod_{\frp} (-p)^{\frac{1}{2}f_\frp
(\phi,\Ind ^G_{I_\frp}u_\frp)}
\] where
here $\frp$ runs over all elements of $S_f(K)$ which divide $p$,
$f_\frp$ is the residue class extension degree of $\frp$ in
$K/\bq$, $I_\frp$ is the inertia subgroup of any fixed place of
$L$ which lies above $\frp$, $u_\frp$ is the augmentation
character of $I_\frp$ and $(-,-)$ is the standard inner product on
$R_G$. \hqed
\end{definition}

We remark that the image of $[O_L,h_{L,\bullet}]$ under $\pi_T^s$
is equal to the element $\chi^s_T(O_L , \det h_\bullet )$ defined
in \cite[\S4.3]{CPT}. This implies that the following result is
equivalent to [loc. cit., Th. 5.9(a)].
\begin{corollary} \label{ecpt}  The class
$\pi^s_T([O_L,h_{L,\bullet}])$ is represented by the homomorphism
which sends each symplectic character $\phi\in \hat G$ to
\[
w_\infty (K,\phi) \mathrm{Pf}(O_L)(\phi) \times
(|G|^{[K:\bq]}|d_K|)^{\frac{\phi(1)}{2}}.
\]
\end{corollary}

\begin{proof} Upon computing the image under $\pi_T^s$ of the class
represented by the explicit homomorphism described in Corollary
\ref{C:kcr}, one finds that it suffices to prove the following result:
if $f$ denotes the homomorphism which sends each symplectic character
$\phi \in \hat G$ to
\[
\lambda_\phi w_\infty(K,\phi)\mathrm{Pf}(O_L)(\phi) \times
|\lambda_\phi \tau(K,\phi)|^{-1},
\]
then $\theta_T^s(f) = 0$.

To prove this we first observe that there exists a finite tamely
ramified abelian extension $N$ of $\bq$ and an element $u$ of $\hat
O_N [G]^\times$ such that, for each symplectic character $\phi \in
\hat G$ one has
\[
\mathrm{Pf}(O_L)(\phi) = \tau(K,\phi)\cdot \Det (u)(\phi).
\]
This can be seen, for example, to be a consequence of the argument
which proves \cite[Proposition 5.11]{CPT}. Hence one has
$\theta_T^s(f) = \theta_T^s(f')$, where $f'$ is the homomorphism which
sends each symplectic character $\phi \in \hat G$ to
\[
\lambda_\phi
 w_\infty(K,\phi)\tau(K,\phi) \times |\lambda_\phi \tau(K,\phi)|^{-1}.
\]
Now the equality \eqref{et} combines with the containment \eqref{pc}
to imply that, for each symplectic character $\phi \in \hat G$, the
real numbers $\lambda_\phi$ and $w_\infty (K,\phi)\tau(K, \phi)$ have
the same sign (note that $\phi = \overline{\phi}$ if $\phi$ is
symplectic). This implies that the element $\Phi$ of $\Hom(R_G^s,(\bq
^c)^{\times})^{\Omega_\bq}$ which sends each symplectic character
$\phi\in \hat G$ to $\lambda_\phi w_\infty(K,\phi)\tau(K,\phi)$
actually belongs to $\Det ^s( \bq [G]^{\times})$.  In addition, since
for each symplectic $\phi\in \hat G$ one has $w_\infty(K,\phi) = \pm
1$, it also implies that $\lambda_\phi w_\infty(K,\phi)\tau(K,\phi) =
|\lambda_\phi \tau(K,\phi)|$.  This in turn implies that $f' =
\Delta^s(\Phi )$, and hence that $\theta_T^s(f ') = 0$, as required.
\end{proof}


\subsection{Torsion modules}\label{Ctors} In this subsection we shall
indicate how certain natural torsion Galois module invariants
introduced by S. Chase (see \cite{Ch}, \cite[Notes to Chapter
III]{Fr0}) are related to the element $\delta^{L}_{O_K[G]}(O_L)$.
 
We first recall that there is a natural $O_L[G]$-equivariant injection
\begin{equation*}
\psi_{L/K}: O_L \otimes_{O_K} O_L \lra
\Map(G,O_L); \,\,\psi_{L/K}(l_1\otimes l_ 2)(g) = g(l_1)l_2.
\end{equation*}

The cokernel $\Cok(\psi_{L/K})$ of $\psi_{L/K}$ is finite, and is of
finite projective dimension as an $O_L[G]$-module because $L/K$ is
tamely ramified. Hence $\Cok(\psi_{L/K})$ determines a class in
$K_0T(O_L[G])$. This class has been studied by Chase and other
authors, and it gives information regarding the $O_K[G]$-module
structure of $O_L$.
 
We write
\[  
i_{L/K}: K_0(O_K [G],L) \lra K_0T(O_L [G])
\]
for the composite of the natural scalar extension morphism
\[  
K_0(O_K [G],L) \lra K_0(O_L [G],L);\quad [X,Y;i] \mapsto [X
\otimes_{O_K} O_L, Y \otimes_{O_K} O_L; i]
\]
and the first isomorphism of (\ref{rt}) with $\A = O_L[G]$.

\begin{lemma}\label{cht} The class of $\Cok(\psi_{L/K})$ in
$K_0T(O_L[G])$ is equal to the image of $\delta_{O_K[G]}^L(O_L)$ under
$i_{L/K}$. 
\end{lemma}

\begin{proof} There is an $O_L[G]$-equivariant
isomorphism
\[ 
\theta: \Map(G,O_L) \stackrel{\sim}\lra
H^L_K\otimes_\Ze O_L; \qquad \theta(\phi) =
(\phi(g))_{g}.
\] 
The stated result follows from the
definition of isomorphism in $\calP_{O_L[G]}\times_L
\calP_{O_L[G]}$ (see \S\ref{fpci}) 
 and the commutativity of the following diagram in $\calP_{L[G]}$:
\begin{equation*}
\begin{CD} 
L\otimes_K L
@> \psi_{L/K} >> \Map (G,L) \\ @V{\id} VV 
@V{\theta}VV\\ L\otimes_K L @>\pi^{L,L}_{K} >>
H_K^L \otimes_{\bz} L
\end{CD}
\end{equation*}
Here $\pi^{L,L}_{K}$ is as defined in Example \ref{class}.
\end{proof}

We remark that Chase has given an explicit description of
$\Cok(\psi_{L/K})$ as an $O_L[G]$-module (see in particular Theorem
2.15 of \cite{Ch}). Hence, Lemma \ref{cht} shows that Chase's
results may be naturally reinterpreted in terms of equivariant
discriminants.
 
\section{Torsors} \label{S:torsors}

We shall now apply the methods developed in this paper to the
`geometrical' setting of Example \ref{geom}.

\subsection{Reduced resolvends} In this subsection, we shall discuss
reduced resolvends coming from torsors of finite group schemes. This
notion was first introduced by L. McCulloh (see \cite{Mc}, and the
references contained therein). Our approach will be rather different
from McCulloh's, and will use ideas that first appeared
in \cite{AP} and \cite{A3}.
\smallskip

We use the notation of Example \ref{geom}. Set $\Gamma:=
G(R^c)=G(F^c)$.

\begin{definition}
We define
\begin{align*}
\bH(\A) &:= \{ \alpha \in \A_{R^c}^{\times}\,|\,
\alpha^{\omega} = g_{\omega} \alpha \,\, \text{for all $\omega \in
\Omega_F$, where $g_{\omega} \in \Gamma$}\}, \\
\cH(\A) &:= \bH(\A) / \Gamma, \\
H(\A) &:= \bH(\A)/(\Gamma \cdot \A^{\times}).
\end{align*}

(Note that we allow the possibility that $R=F$, in which case we have
$\A=A$.) \hqed
\end{definition}
 
We suppose for the moment that $R$ is a local ring, and we let $\pi: X
\to \Spec(R)$ be a $G$-torsor. Then the line bundle $\L_{\pi}$
associated to $\pi$ (see Example \ref{geom}) is isomorphic to the
trivial line bundle on $G^*$. Let $s_{\pi}: \A \xrightarrow{\sim}
\L_{\pi}$ be any trivialisation of $\L_{\pi}$. Recall from Example
\ref{geom} that over $\Spec(R^c)$, $\pi$ becomes isomorphic to the
trivial torsor $\pi_0$, i.e. there is an isomorphism $i: X \otimes_R
R^c \xrightarrow{\sim} G \otimes_R R^c$ of schemes with
$G$-action. This isomorphism is not unique. If $i'$ is any other such
isomorphism, then $i^{-1}i': G \otimes_R R^c \xrightarrow{\sim} G
\otimes_R R^c$ is also an isomorphism of schemes with $G$-action, and
so is given by translation by an element of $\Gamma$.

Let $\xi_\pi: \L_\pi \otimes_R R^c \xrightarrow{\sim} \fA \otimes_R
R^c$ be the splitting isomorphism corresponding to $i$, and consider
the map from $\A_{R^c}$ to itself defined by
\begin{equation*}
\A_{R^c} \xrightarrow{s_{\pi} \otimes_R R^c} \L_{\pi}
\otimes_R R^c \xrightarrow{\xi_{\pi}} \A_{R^c}.
\end{equation*}
This is an isomorphism of $\A_{R^c}$-modules, and so it is given by
multiplication by some element $\r(s_{\pi})$ in
$\A_{R^c}^{\times}$. We refer to $\r(s_{\pi})$ as a {\it resolvend} of
$s_{\pi}$. Note that $\r(s_{\pi})$ depends upon $s_{\pi}$ as well as
upon the choice of $\xi_{\pi}$.

If $\omega \in \Omega_F$, then $\xi_{\pi}^{\omega} =
g_{\omega}\xi_{\pi}$, for some $g_{\omega} \in \Gamma$. Since
$s_{\pi}^{\omega} = s_{\pi}$, we deduce that $\r(s_{\pi})^{\omega} =
g_{\omega}\r(s_{\pi})$, and so $\r(s_{\pi}) \in \bH(\A)$. As $i$ is
only well-defined up to translation by an element of $\Gamma$, it
follows that changing our choice of $\xi_{\pi}$ alters $\r(s_{\pi})$
by multiplication by an element of $\Gamma$ (see also \cite[Lemmas 1.2
and 1.3]{BT}). Hence the image $r(s_{\pi})$ of $\r(s_{\pi})$ in
$\cH(\A)$ depends only upon the trivialisation $s_{\pi}$. We refer to
$r(s_{\pi})$ as the {\it reduced resolvend} of $s_{\pi}$.

We next observe that changing our choice of trivialisation $s_{\pi}$
alters $\r(s_{\pi})$ by multiplication by an element of
$\A^{\times}$. It follows that the image of $r(s_{\pi})$ in the group
$H(\A)$ depends only upon the isomorphism class of the $G$-torsor
$\pi$. 

The following result was first proved by McCulloh in the case in which
$G$ is a constant group scheme \cite[\S1 and \S2]{Mc}. In the
generality given below, it is due to Byott \cite[Lemma 1.11,
Proposition 2.12]{By}. A different proof of this result is given in
\cite{A3}.

\begin{theorem} \label{T:byott}
Suppose that $R$ is a local ring. Then the map
\begin{equation*}
H^1(\Spec(R), G) \to H(\A);\qquad \pi
\mapsto [r(s_{\pi})]
\end{equation*}
is an isomorphism. \hfill$\square$
\end{theorem}
\smallskip

We shall now show that, for any field $F$, the subgroup
$H(A)$ of $A_{F^c}^{\times}/(\Gamma \cdot
A^{\times})$ has a natural functorial description. In order to do
this, we first make the following definition.

\begin{definition} Let $R$ be any commutative ring, and suppose that
$\F$ is any contravariant functor from the category of finite,
flat, commutative group schemes over $\Spec(R)$ to the category of
abelian groups.

If $G$ is any finite, flat, commutative group scheme over $R$, let
$m:G \times G \to G$ denote multiplication on $G$, and write $p_i:G
\times G \to G$ ($i=1,2$) for projection onto the $i$-th factor. Let
$m^*,p_{i}^{*}: \F(G) \to \F(G \times G)$ denote the homomorphisms
induced by $m$ and $p_i$ respectively. An element $x \in \F(G)$ is
said to be {\it primitive} if $m^*(x) = p_{1}^{*}(x) \cdot
p_{2}^{*}(x)$. (Note in particular that it follows that the subset of
primitive elements of $\F(G)$ is a group.)

It is sometimes helpful to formulate this definition in terms of Hopf
algebras rather than group schemes. This may be done as follows. Each
functor $\F$ as above naturally corresponds to a covariant functor
$\overline{\F}$ from the category of finite, flat, commutative and
cocommutative Hopf algebras over $R$ to the category of abelian
groups. For any such Hopf algebra $\cA$, let $\Delta: \cA \to \cA
\otimes_R \cA$ denote the comultiplication on $\cA$. Define $i_1,i_2:
\cA \to \cA \otimes_R \cA$ by $i_1(a) = a \otimes 1$ and $i_2(a) = 1
\otimes a$ for each $a \in \cA$. Write $\Delta^*, i_{1}^{*},
i_{2}^{*}: \overline{\F}(\cA) \to \overline{\F}(\cA \otimes_R
\cA)$ for the homomorphisms induced by $\Delta$, $i_1$, and $i_2$
respectively. Then an element $y \in \overline{\F}(\cA)$ is
primitive if $\Delta^*(y) = i_{1}^{*}(y) \cdot i_{2}^{*}(y)$. \hqed
\end{definition}

\begin{theorem} \label{T:primres}
The group $H(A)$ is equal to the subgroup of primitive elements of the
group $A_{F^c}^{\times}/(\Gamma \cdot A^{\times})$.
\end{theorem}

\begin{proof} Suppose first that $a \in \bH(A)$, with $a^{\omega} =
g_{\omega}a$ for all $\omega \in \Omega_F$. We have that $A_{F^c} =
F^c[\Gamma]$, and the comultiplication map $\Delta: A_{F^c} \to A_{F^c}
\otimes_F A_{F^c}$ is induced by $\Delta(\gamma) = \gamma \otimes
\gamma$ for $\gamma \in \Gamma$. This implies that, for all $\omega
\in \Omega$,
\begin{equation*}
\Delta(a)^{\omega} = (g_{\omega} \otimes g_{\omega}) \Delta(a),\quad
i_1(a)^{\omega} = (g_{\omega} \otimes 1) i_1(a),\quad
i_2(a)^{\omega} = (1 \otimes g_{\omega}) i_2(a).
\end{equation*}
Hence we have that $\Delta(a)[i_1(a)i_2(a)]^{-1} \in (A_{F^c}
\otimes_F A_{F^c})^{\Omega_F} = A \otimes_F A$. This implies that the
class of $a$ in $A_{F^c}^{\times}/(\Gamma \cdot A^{\times})$ is
primitive.

Suppose conversely that $a \in A_{F^c}^{\times}$ represents a
primitive class in $A_{F^c}^{\times}/(\Gamma \cdot A^{\times})$. For
each $\omega \in \Omega_F$, write $a^{\omega} = u_{\omega}a$. We wish
to show that $u_{\omega} \in \Gamma$ for all $\omega$.

Since $a \in A_{F^c}^{\times}$ represents a primitive
class in $A_{F^c}^{\times}/(\Gamma \cdot A^{\times})$,
we have
\begin{equation*}
\frac{\Delta(a)}{i_1(a)i_2(a)} = (g \otimes h) \beta,
\end{equation*}
where $g,h \in \Gamma$, and $\beta \in (A \otimes_F
A)^{\times}$. Hence if $\omega \in \Omega_F$, then $\beta^{\omega} =
\beta$ and
\begin{equation} \label{E:prim 1}
\left[ \frac{\Delta(a)}{i_1(a)i_2(a)} \right]^{\omega}
=(g^{\omega} \otimes h^{\omega}) \beta = (g^{\omega}g^{-1} \otimes
h^{\omega}h^{-1}) \left(\frac{\Delta(a)}{i_1(a)i_2(a)}\right).
\end{equation}

On the other hand, we also have that
\begin{equation} \label{E:prim 2}
\left[ \frac{\Delta(a)}{i_1(a)i_2(a)} \right]^{\omega}
=
\frac{\Delta(u_{\omega})}{i_1(u_{\omega})i_2(u_{\omega})} 
\left(\frac{\Delta(a)}{i_1(a)i_2(a)}\right).
\end{equation}
Hence \eqref{E:prim 1} and \eqref{E:prim 2} imply that
\begin{equation}
\frac{\Delta(u_{\omega})}{i_1(u_{\omega})i_2(u_{\omega})} =
g^{\omega}g^{-1} \otimes h^{\omega}h^{-1}.
\end{equation}

Therefore, in order to prove that $u_{\omega} \in \Gamma$, it suffices
to prove the following assertion. Suppose that $u \in
A_{F^c}^{\times}$ satisfies
\begin{equation} \label{E:primu}
\frac{\Delta(u)}{i_1(u) i_2(u)} = e \otimes f,
\end{equation}
with $e,f \in \Gamma$. Then $u \in \Gamma$.

In order to establish this, we argue as follows. Suppose first that
$u=a_{\gamma}\gamma$ for some $\gamma \in \Gamma$. Then
\begin{equation*}
\frac{\Delta(u)}{i_1(u)i_2(u)} = \frac{a_{\gamma}(\gamma \otimes
\gamma)}{a_{\gamma}^{2}(\gamma \otimes \gamma)} = \frac{1}{a_{\gamma}}
1_{\Gamma} \otimes 1_{\Gamma},
\end{equation*}
and so in this case $u$ satisfies \eqref{E:primu} if and only if
$a_{\gamma} = 1$. Hence we may assume that $u = \sum_{\gamma \in
\Gamma} a_{\gamma} \gamma$ with $a_{\gamma_1} a_{\gamma_2} \neq 0$ for
some $\gamma_1 \neq \gamma_2$.

Under this assumption, \eqref{E:primu} implies that
\begin{align}
\sum_{\gamma} a_{\gamma} \gamma \otimes \gamma &= (e \otimes f)
\sum_{\beta, \beta'} a_{\beta} a_{\beta'} (\beta \otimes \beta') \notag \\ 
&=
(e \otimes f) [a_{\gamma_1}^{2}(\gamma_1 \otimes \gamma_1) +
a_{\gamma_2}^{2} (\gamma_2 \otimes \gamma_2) + a_{\gamma_1}
a_{\gamma_2} (\gamma_1 \otimes \gamma_2) + (\textnormal{other
terms})]. \label{E:groupeq}
\end{align}
If $e=f$, then \eqref{E:groupeq} implies that $a_{\gamma_1}
a_{\gamma_2} (\gamma_1 \otimes \gamma_2)=0$, whence it follows that
$a_{\gamma_1} a_{\gamma_2}=0$, which is a contradiction. On the other
hand, if $e \neq f$, then we see from \eqref{E:groupeq} that
$a_{\gamma_1}^{2} = a_{\gamma_2}^{2} =0$, which is also a
contradiction. Hence $u \in \Gamma$ as claimed.

This completes the proof of the result.
\end{proof}

\subsection{Class invariant maps}

We retain the notation established in the previous section.
In this subsection we describe a natural refinement of the class
invariant homomorphism introduced by Waterhouse in \cite{W}.

Suppose that $F$ is a number field. Then the class invariant
homomorphism $\psi$ is defined by
\begin{equation*}
\psi: H^1(\Spec(R),G) \to \Pic(G^*);\qquad \pi \mapsto (\L_{\pi}).
\end{equation*}

\begin{remark} \label{R:psi}
In general $\psi$ is not injective (see e.g. \cite{A4}, \cite{P1},
\cite{ST}). It is known that the image of $\psi$ is always contained
in the subgroup $\PPic(G^*)$ of primitive classes of $\Pic(G^*)$ (see
e.g. \cite{CM}).  It is shown by the first-named author in \cite{A2}
that if the composition series of every geometric fibre of $G \to
\Spec(R)$ of residue characteristic $2$ does not contain a factor of
local-local type, then the image of $\psi$ is in fact equal to
$\PPic(G^*)$.  \hqed
\end{remark}

Now suppose that $\pi:X \to \Spec(R)$ is a $G$-torsor, and fix a
choice of isomorphism $\xi_{\pi}$ as in \eqref{E:splitting}.  For ease
of notation, we use the same symbol $\xi_{\pi}$ for the map $\L_{\pi}
\otimes_R F^c \xrightarrow{\sim} A_{F^c}$ induced by
\eqref{E:splitting}. In general the element $[\L_{\pi}, \A; \xi_{\pi}]
\in K_0(\A, R^c)$ depends upon the choice of $\xi_{\pi}$.

Recall that there is a natural morphism
\begin{equation*}
\iota_{\A,R^c}: K_0(\A,R^c) \rightarrow K_0(\A,F^c)
\end{equation*}
arising via extension of scalars from $R$ to $F$. The image of
$\iota_{\A,R^c}([\calL_\pi,\A;\xi_\pi])$ under the isomorphism
\eqref{E:idelic} with $E = F^c$ may be described as follows. Choose
any trivialisation $s_{\pi}: A \xrightarrow{\sim} \L_{\pi} \otimes_R
F$. For each finite place $v$ of $F$, choose a trivialisation
$t_{\pi,v}: \A_v \xrightarrow{\sim} \L_{\pi} \otimes_R R_v$. It is
easy to check that the isomorphism
\begin{equation*}
A_v \xrightarrow{t_{\pi,v}} \L_{\pi} \otimes_R F_v
\xrightarrow{(s_{\pi} \otimes_F F_v)^{-1}} A_v
\end{equation*}
is given by multiplication by $\r(t_{\pi,v}).\r(s_{\pi})^{-1} \in
A_{v}^{\times}$. Then Remark \ref{exdes} implies that a representative
of the image of $\iota_{\A,R^c}([\L_{\pi},\A ;\xi_{\pi}])$ under the
isomorphism \eqref{E:idelic} is given by
\begin{equation} \label{E:representative}
\prod_{v\in S_f(F)}(\r(t_{\pi,v}).\r(s_{\pi})^{-1}) \times
\r(s_{\pi}) \in J_f(A) \times A_{F^c}^{\times}.
\end{equation}

Observe that if we change the choice of splitting isomorphism
$\xi_{\pi}$, the element $\r(s_{\pi})$ is altered by
multiplication by an element of $\Gamma$. This motivates the
following definition:

\begin{definition}
Let $\Delta''_{\A,F^c}$ denote the composition of the homomorphism
\begin{equation*}
\Delta'_{\A,F^c}: A^{\times} \longrightarrow \frac{J_f(A)}{U_f(\A)} \times
A^{\times}_{F^c}
\end{equation*}
defined in Example \ref{exdes1}(ii) with the quotient map
\begin{equation*}
\frac{J_f(A)}{U_f(A)} \times A^{\times}_{F^c} \longrightarrow
\frac{J_f(A)}{U_f(A)} \times \frac{A^{\times}_{F^c}}{\Gamma}.
\end{equation*}
We define $\oK_0(\A,F^c)$ to be the quotient of $K_0(\A,F^c)$ defined by
the isomorphism
\begin{equation} \label{E:idelic K}
\oK_0(\A,F^c) \simeq \frac{\dfrac{J_f(A)}{U_f(\A)} \times
\dfrac{A_{F^c}^{\times}}{\Gamma}}{\Image(\Delta'' _{\A,F^c})}. 
\end{equation}
We write $\oK_0(\A,R^c)$ for the image of $K_0(\A,R^c)$ under the
composite of the scalar extension map $\iota_{\A,R^c}$ and the natural
projection map $K_0(\A,F^c) \to \oK_0(\A,F^c)$.  \hqed
\end{definition}

\begin{lemma}\label{oKl} The natural projection map
 $K_0(\mathfrak A, F^c) \rightarrow \oK_0(\mathfrak A, F^c)$ is
bijective if and only if the points of $\Gamma = G(F^c)$ are fixed by
$\Omega_F$.
\end{lemma}

\begin{proof} It is clear that the stated projection is bijective if
and only if $\Gamma \subseteq {\mathfrak A}^{\times}_v$ for all $v \in
S_f(F)$, i.e. if and only if $\Gamma \subseteq {\mathfrak
A}^{\times}$. But this last condition obtains if and only if the
points of $\Gamma$ are fixed by $\Omega_F$.
\end{proof}

It follows from \eqref{E:representative} that the class of
$\iota_{\A,R^c}([\L_{\pi},\A ;\xi_{\pi}])$ in $\oK_0(\A,F^c)$ is
represented by
\begin{equation*}
\prod_{v\in S_f(F)}(\r(t_{\pi,v}).\r(s_{\pi})^{-1}) \times
\r(s_{\pi}) \in J_f(A) \times \cH(A).
\end{equation*}
It is easy to check that this class is independent of the choice
of splitting isomorphism $\xi_{\pi}$ and that it depends only upon the
isomorphism class of the $G$-torsor $\pi: X \to \Spec(R)$. 

\begin{theorem} \label{T:kappa}
The map
\begin{equation*}
\opsi: H^1(\Spec(R),G) \to \oK_0(\A,R^c);\qquad \pi \mapsto
\iota_{\A,R^c}([\L_{\pi},\A;\xi_{\pi}])
\end{equation*}
is an injective group homomorphism.
\end{theorem}

\begin{proof} Let us first show that $\opsi$ is a group
homomorphism. Suppose that $\pi$ and $\pi'$ are $G$-torsors. Let
$\xi_{\pi}$ and $\xi_{\pi'}$ denote splitting isomorphisms associated
to $\pi$ and $\pi'$ respectively. It follows from the functoriality of
Waterhouse's construction in \cite{W} that there is a canonical
isomorphism $\L_{\pi} \otimes_{\A} \L_{\pi'} \simeq \L_{\pi \cdot
\pi'}$ and that the map 
\begin{equation*}
\xi_\pi\otimes_{\A_{R^c}}\xi_{\pi '}: (\calL_\pi \otimes_R R^c)
 \otimes_{\A_{R^c}}(\calL_{\pi '}\otimes_R R^c)
 \stackrel{\sim}\longrightarrow \A_{R^c}\otimes_{\A_{R^c}}\A_{R^c}
 \stackrel{\sim}\longrightarrow \A_{R^c}
\end{equation*}
is a splitting isomorphism associated to $\pi \cdot \pi'$. Since
in $K_0(\A,R^c)$, we have
\begin{equation*}
[\L_{\pi},\A ;\xi_{\pi}] + [\L_{\pi'},\A ;\xi_{\pi'}] = [\L_{\pi}
\otimes_{\A} \L_{\pi'}, \A ; \xi_{\pi} \otimes_{\A_{R^c}} \xi_{\pi'}],
\end{equation*}
(cf. for example \cite[proof of Lemma 2.6(i)]{BB1}), it follows that
\begin{equation*}
\opsi(\pi) + \opsi(\pi') = \opsi(\pi \cdot \pi'),
\end{equation*}
as required.

To show that $\opsi$ is injective, we argue as follows. If we compose
the natural injection \break $H^1(\Spec(R),G) \to H^1(F,G)$ with the
isomorphism $H^1(F,G) \to H(A)$ afforded by Theorem
\ref{T:byott}, then we obtain an injection $H^1(\Spec(R),G) \to
H(A)$. It follows easily from the definitions that this
injection is the same as the homomorphism obtained by composing
$\opsi$ with the inclusion $\oK_0(\A,R^c) \subseteq \oK_0(\A,F^c)$ and the
natural projection
\begin{equation*}
\oK_0(\A,F^c) \to \frac{A_{F^c}^{\times}}{\Gamma \cdot A^{\times}}.
\end{equation*}
This implies that $\opsi$ is injective.
\end{proof}
\smallskip

\begin{remark}\label{comp} Let $\phi \in \hat\Gamma$. Then for each
$\gamma \in \Gamma$ and $v \in S_\infty (F)$, the complex number
 $\sigma_v(\phi (\gamma ))$ is a root of unity and hence has absolute
 value $1$. This implies that the morphism 
\begin{equation*}
 \partial_{\A, F^c} : K_0(\A,F^c) \to \mathrm{AC}(\A)
\end{equation*}
defined in \S\ref{S:concat} induces a morphism $\oK_0(\A,F^c) \rightarrow
 \mathrm{AC}(\A)$. Upon composing this morphism with $\opsi$ we obtain
a morphism
\begin{equation*}
\hat{\psi}:H^1(\Spec(R),G) \to \mathrm{AC}(\A) \simeq \hPic(G^*)
\end{equation*}
which was first introduced and studied by the first-named author and
G. Pappas in \cite{AP}. In turn, upon composing $\hat{\psi}$ with the
natural surjection $\mathrm{AC}(\A) \to \Pic(G^*)$ we obtain the
homomorphism
\begin{equation*}
\psi :H^1(\Spec(R),G) \to \Pic(G^*)
\end{equation*}
which was introduced by Waterhouse in \cite{W}. \hqed
\end{remark}

Write $P\oK_0(\A,R^c)$ for the subgroup of primitive elements of
$\oK_0(\A,R^c)$.

\begin{theorem} \label{T:primclass}
(i) The image of $\opsi$ lies in $P\oK_0(\A,R^c)$.

(ii) Suppose that the composition series of every geometric fibre of
$G \to \Spec(R)$ of residue characteristic $2$ does not contain a
factor of local-local type. Then the image of $\opsi$ is equal to
$P\oK_0(\A,R^c)$. 
\end{theorem}

\begin{proof}
(i) Recall from Example \ref{geom} that we have a canonical
isomorphism 
\begin{equation*}
\Ext^1(G^*,\bG_m) \simeq H^1(\Spec(R), G).
\end{equation*} 
Composing this isomorphism with $\opsi$ gives a homomorphism
\begin{equation*}
\vk_{G^*}: \Ext^1(G^*,\bG_m) \to \oK_0(\A,R^c)
\end{equation*}
whose image is the same as that of $\opsi$. It follows from the
functoriality of Waterhouse's construction in \cite{W} that the
homomorphism $\vk_{G^*}$ is functorial in $G^*$.

Let $m: G^* \times G^* \to G^*$ denote multiplication on $G^*$,
and write $p_i:G^* \times G^* \to G^*$ ($i=1,2$) for projection
onto the $i$-th factor. Let
\begin{equation*}
m^*,p_{i}^{*}: \Ext^1(G^*,\bG_m) \to \Ext^1(G^* \times G^* ,\bG_m)
\end{equation*}
be the homomorphisms induced by $m$ and $p_i$. Then it follows from
standard properties of the functor $\Ext^1$ that $m^*(x) =
p_{1}^{*}(x) + p_{2}^{*}(x)$ for all $x \in \Ext^1(G^*,\bG_m)$
(see e.g. \cite[Proposition 1 on p. 163, and p. 182]{Serre}). Since
$\vk_{G^*}$ is functorial in $G^*$, this implies that
$\vk_{G^*}(x)$ is a primitive element of $\oK_0(\A,R^c)$.
\smallskip

(ii) It is shown in \cite{A3} that there is a natural isomorphism
$\Ker(\psi) \simeq H(\A)$. From \cite[Theorem 1.3]{A2}, we know that,
under the given hypotheses, we have that $\Im(\psi) =
\PPic(G^*)$. Since $\opsi$ is injective, this implies that there is an
exact sequence
\begin{equation} \label{E:image}
0 \to H(\A) \to \Im(\opsi)
\to  \PPic(G^*) \to 0.
\end{equation}

Next, we observe that upon comparing the exact sequences \eqref{krel}
with $\Lambda = R^c$ and $\Lambda = F^c$, we obtain a commutative
diagram of short exact sequences
\begin{equation*}
\begin{CD} 
0 @> >> \cok[K_1(\A) \rightarrow K_1(\A_{R^c})] @> >>
K_0(\A,R^c) @> >> \ker[K_0(\A) \rightarrow K_0(\A_{R^c})]@> >>
0\cr @. @V VV @V \iota_{\A,R^c} VV @V \subseteq VV \cr
 0@> >> \cok[K_1(\A) \rightarrow K_1(A_{F^c})] @> >> K_0(\A,F^c) @> >>
 \ker[K_0(\A) \rightarrow K_0(A_{F^c})]@> >> 0.
\end{CD}
\end{equation*} 

Now, since $\A$ is commutative, the second and fourth terms of the
lower sequence can be identified with $A_{F^c}^\times/\A^\times$
and $\Pic(\A)$ respectively (cf. the end of Remark \ref{ct}). Further,
 the image of the natural composite morphism
\begin{equation*}
K_1(\A_{R^c}) \rightarrow K_1(A_{F^c}) \rightarrow \cok[K_1(\A)
\rightarrow K_1(A_{F^c})] \cong A_{F^c}^\times/\A^\times
\end{equation*}
is equal to $\A_{R^c}^\times/\A^\times$. Upon identifying $\Pic(\A)$
with $\Pic(G^*)$ the above diagram implies that there is an exact
sequence
\[ 
1 \rightarrow \frac{\A^{\times}_{R^c}}{\A^\times} \rightarrow
 \iota_{\A,R^c}(K_0(\A,R^c)) \rightarrow \Pic(G^*).
\]
It is easy to see that this in turn induces an exact sequence
\begin{equation} \label{E:barimage}
1 \to \frac{\A_{R^c}^{\times}}{\Gamma \cdot \A^{\times}} \to \oK_0(\A,
R^c) \to \Pic(G^*).
\end{equation}

Since $\Im(\opsi) \subseteq P\oK_0(\A,R^c)$ it follows from
\eqref{E:image} that the natural map $P\oK_0(\A,R^c) \to \PPic(G^*)$
is surjective. Hence Theorem \ref{T:primres} implies that
\eqref{E:barimage} induces an exact sequence
\begin{equation} \label{E:primexact}
0 \to H(\A) \to P\oK_0(\A,R^c)
\to  \PPic(G^*) \to 0.
\end{equation}
We now see from \eqref{E:image} and \eqref{E:primexact} that
$\Im(\opsi) = P\oK_0(\A,R^c)$, as asserted.
\end{proof}

\subsection{Rings of integers}

In this subsection we shall briefly explain how McCulloh's results on
realisable classes of tame extensions in \cite{Mc} can be naturally
lifted to the corresponding relative algebraic $K$-group. Proposition
\ref{cht} implies that such realisability results in the relative
algebraic $K$-group are related to those of Chase in \cite{Ch}.
\smallskip

Suppose that $F$ is a number field with ring of integers $R$. Let $G$
be any finite, constant abelian group scheme over $\Spec(R)$, and set
$\Gamma := G(F^c)$. Then in this case, we have $G = \Spec(\B)$ and
$G^* = \Spec(\A)$, where
\begin{equation*}
\B = \Map(\Gamma, R),\qquad \A = R[\Gamma].
\end{equation*}
Write $G_F$ and $G_F^*$ for the generic fibres of $G$ and $G^*$
respectively. Note that, since $G_F$ is a constant group scheme over
$\Spec(F)$, Lemma \ref{oKl} implies that $\oK_0(\A,F^c)$ is
canonically isomorphic to $K_0(\A,F^c)$.

Now suppose that $\pi:X \to \Spec(F)$ is a $G_F$-torsor, and write
$L_{\pi}$ for its associated line bundle on $G^*_F$. Then
$X=\Spec(C_{\pi})$, where $C_{\pi}$ is Galois algebra extension of
$F$ with Galois group $\Gamma$. 

\begin{definition} \label{D:tame}
We say that $\pi$ is {\it tame} if it is trivialised by a tamely
ramified extension of $F$.  Equivalently, $\pi$ is tame if and only if
the extension $C_{\pi}/F$ is tamely ramified. 

We write $H^1_t(\Spec(F),G_F)$ for the subgroup of $H^1(\Spec(F),G_F)$
consisting of isomorphism classes of tame $G_F$-torsors. \hqed
\end{definition}

Let $O_{\pi}$ denote the integral closure of $R$ in $C_{\pi}$.  If
$\pi$ is tame, then Noether's theorem (see e.g. \cite[Chapter I,
\S3]{Fr0}) implies that $O_{\pi}$ is a locally free, rank one
$\A$-module, and so it yields a line bundle on $G^*$ which we shall
denote by $\cM_{\pi}$. Let $\pi_0: G_F \to \Spec(F)$ denote the
trivial $G_F$-torsor, and set
\begin{equation*}
\oL_{\pi} := \cM_{\pi} \otimes \cM_{\pi_0}^{-1}.
\end{equation*}

It follows from the definitions that there is a natural
isomorphism $\oL_{\pi} \otimes_{R} F \simeq L_{\pi}$. If we make
this identification, then any choice of splitting isomorphism
$\xi_{\pi}$ for $\pi$ gives a class $[\oL_{\pi},\A ;\xi_{\pi}]$ in
$K_0(\A,F^c)$. This class is
independent of the choice of $\xi_{\pi}$. Hence we obtain a map
\begin{equation}
{\opsi_t}: H^1_t(\Spec(F),G_F) \to K_0(\A,F^c); \qquad \pi \mapsto
[\oL_{\pi},\A ;\xi_{\pi}].
\end{equation}

If $F$ is a number field, then a representative in $J_f(A) \times
A_{F^c}^{\times}$ of ${\opsi_t}(\pi)$ may be given exactly as in the
previous section. It may be shown via an argument very similar to that
given in Theorem \ref{T:kappa} that the map ${\opsi_t}$ is
injective. However, in general, ${\opsi_t}$ is not a group
homomorphism (see Remark \ref{R:misc} below).

\begin{definition} \label{D:realisable}
We say that an element of $K_0(\A,F^c)$ is {\it realisable} if it lies
in the image of ${\opsi_t}$. \hqed
\end{definition}

The general strategy for describing the set of realisable classes in
$K_0(\A,F^c)$ is as follows. Let
\begin{equation*}
\bj: J_f(A) \times [A_{F^c}^{\times}/\Gamma] \to
K_0(\A,F^c)
\end{equation*}
denote the obvious quotient map afforded by \eqref{E:idelic K}
(together with Lemma \ref{oKl}), and write $\cH(\bbA(A))$ for the
restricted direct product of the groups $\cH(A_v)$ with respect to the
subgroups $\cH(\A_v)$. Then the obvious natural maps $A_{v}^{\times}
\to \cH(A_v)$ induce a homomorphism
\begin{equation*}
\rag: J_f(A) \to \cH(\bbA(A)).
\end{equation*}
Suppose that ${\underline c}:= c_{\fin} \times c_{\infty} \ \in J_f(A)
\times [A_{F^c}^{\times}/\Gamma]$ satisfies $\bj({\underline c}) \in
\Im({\opsi_t})$. McCulloh shows that the idele $\rag(c_{\fin})$ admits
a special type of decomposition (see \cite[Theorem 5.6]{Mc}). It is a
straightforward exercise to use the results of \cite{Mc} to
characterise the image of $\opsi_t$ in a manner similar to that given
in \cite[Theorem 6.7]{Mc}. We omit the details in order to keep this
paper to a reasonable length.

\begin{remarks} \label{R:misc}
(i) McCulloh has shown that the image of the map obtained by composing
$\opsi_t$ with the quotient map $K_0(\fA,F^c) \to \Pic(G^*)$ is a
subgroup of $\Pic(G^*)$ (see \cite[Corollary 6.20]{Mc}). On the other
hand, by using the description of the image of $\opsi_t$ mentioned
above, it may be shown that in general, $\Im({\opsi_t})$ is not a
group. This situation is very similar to that which obtains in the
case of realisable classes in classgroups of sheaves discussed in
\cite[Remark 2.10(iii)]{AB1}.

(ii) By composing ${\opsi_t}$ with the natural morphism
$\partial_{\A,F^c}: K_0(\A,F^c) \rightarrow \mathrm{AC}(\A)$ defined in
\S\ref{S:concat}, we obtain a map
\begin{equation*}
\hat{\psi_t}: H^{1}_{t}(\Spec(F),G_F) \to \mathrm{AC}(\A).
\end{equation*}
The results of \cite{Mc} may also be used to give a description of
$\Im(\hat{\psi_t})$. It follows from Proposition \ref{pr} that such a
description may be viewed as a realisability result for certain
metrised $\A$-modules (cf. Remark \ref{comp}).  \hqed
\end{remarks}

\end{document}